\documentclass[11 pt]{article}
\usepackage{amsfonts,amsmath,amssymb,amsthm,enumerate,color,hyperref,bbm,tabularx,ifthen,twoopt,mathabx}
\usepackage[T1]{fontenc}
\usepackage[francais,english]{babel}
\usepackage{graphicx}
\usepackage[utf8]{inputenc}
\usepackage{subfig}
\usepackage{array}
\usepackage{stmaryrd}
\usepackage{enumitem}
\usepackage{natbib}
\usepackage{ulem}

\newcounter{hypA}

\def\PE{\mathbb{E}}

\def\Var{\mathop{\rm Var}\nolimits}

\def \1{\mathbbm{1}}
\def\Y{\mathbf{Y}}
\def\betabold{\boldsymbol{\beta}}
\def\X{\mathbf{X}}
\def\u{\mathbf{u}}
\def\Z{\mathbf{Z}}
\def\R{\mathbf{R}}
\def\W{\mathbf{W}}

\def\G{\mathbf{G}}

\def\g{\mathbf{g}}

\def\e{\mathbf{e}}

\def\l{\mathbf{l}}

\def\Id{\textrm{Id}}
\def\rset{\mathbb{R}}

\newtheorem{theorem}{Theorem}
\newtheorem{prop}{Proposition}
\newtheorem{lemma}{Lemma}

\theoremstyle{remark}

\newtheorem{assumption}{Assumption}

\begin{document}


\title{ Heritability estimation of diseases in case-control studies}

\author{Anna Bonnet}

\maketitle

\begin{abstract}
In the field of genetics, the concept of heritability refers to the proportion of variations of a biological trait or disease that can be explained by genetic factors. 
Quantifying the heritability of a disease is a fundamental challenge in human genetics, especially when the causes are plural and not clearly identified. Although the literature regarding heritability estimation for binary traits is less rich than for quantitative traits, several methods have been proposed to estimate the heritability of complex diseases. However, to the best of our knowledge, the existing methods are not supported by theoretical grounds. Moreover, most of the methodologies do not take into account a major specificity of the data coming from medical studies, which is the oversampling of the number of patients compared to controls. We propose in this paper to investigate the theoretical properties of the method developed by \cite{golan2014}, which is very efficient in practice, despite the oversampling of patients. Our main result is the proof of the consistency of this estimator. We also provide a numerical study to compare two approximations leading to two heritability estimators.
\end{abstract}

\section{Introduction}

In the field of genetics, the concept of heritability refers to the proportion of variations of a biological trait or disease that can be explained by genetic factors.  Quantifying the heritability is a major challenge for diseases that are suspected to have a strong genetic component but which causes are often vague and multiple.
Indeed, determining a high value of heritability is a powerful argument in favor of further research for genetic causes, but it also opens the possibility of predicting a risk of illness based on the genetic background.
\\
There exist several methods to estimate the heritability of quantitative traits, which we will describe hereafter, with interesting theoretical and practical properties. Regarding binary traits, such as the presence or absence of a disease, a few methodologies have been proposed, but as far as we know, none of them has been validated theoretically. \cite{golan2014} developed a method that they compared to recent methodologies and which was shown to be very efficient in practice. The aim of this paper is to investigate the theoretical properties of \cite{golan2014}'s method.
\\
Let us first recall the main existing methods to estimate the heritability of quantitative traits, which will be strongly linked to the methods used for binary traits.
Linear Mixed Models (LMMs) have been widely used for estimating the heritability of quantitative traits. Indeed,  \citet{Yang_2010} proposed for instance to estimate the heritability of human height by using a classical LMM defined by 

\begin{equation}
\Y= \X\betabold + \Z\u +\e,
\label{eq:LMM}
\end{equation}
where $\Y=(\Y_1,\dots,\Y_n)'$ is the vector of observations of a phenotype of interest, $\X$ is a $n\times p$ matrix of predictors (or fixed effects), $\betabold$
is a $p\times 1$ vector containing the unknown linear effects of the predictors, and $\u$ and $\e$ correspond respectively to the genetic and the environmental random effects. We assume that $\u$ and $\e$ are Gaussian random
effects with variances $\sigma_u^{\star 2}$ and $\sigma_e^{\star 2}$ respectively. Moreover, $\Z$ is a $n\times N$
matrix which contains the genetic information.
They proposed to estimate the parameter 

\begin{equation}
\eta^\star=\frac{N\sigma_u^{\star 2}}{N\sigma_u^{\star 2}+\sigma_e^{\star2}},
\label{eq:eta_star_intro}
\end{equation}
commonly considered as the mathematical definition for heritability since it determines how 
the variance is shared between $\u$ and $\e$.
\\
Several methods were developed to estimate the parameter $\eta^\star$, see
\cite{patterson_thompson_1971},  \cite{searle_casella_mcculloch_1992},   \cite{yang_lee_goddard_visscher_2011},  \citet{pirinen_donnelly_spencer}, \citet{zhou2012}.
\\
From a theoretical point of view, \citet{nous} showed the asymptotic normality of the maximum likelihood estimator of $\eta^\star$ as well as a central limit theorem leading to confidence intervals for $\eta^\star$.
\\
\bigskip
\\
The previous modeling and the corresponding methods obviously do not apply when considering non continuous traits. However, the quantitative and the binary cases can be related by assuming the existence of an underlying Gaussian variable linked to the binary phenotype. More precisely, there are two main modelings which connect binary phenotypes to a continuous and unobserved quantity called the liability. 
\\
The first one consists in assuming that the observations $\Y_1, \dots, \Y_n$ are distributed according to the following Generalized Linear Mixed Model (GLMM): 
\begin{equation}
\Y_i \sim \mathcal{B}(p_i)
\label{eq:GLMM}
\end{equation}
with  $p_i=g(\l_i)$ where g is a link function and $\l_i$ is defined as

\begin{equation}
\label{eq:liability}
\l=\Z \u + \e
\end{equation}
\\
with $\l=(\l_1,\dots, \l_n)$, $\u \sim \mathcal{N}(0,\sigma_u^{\star2})$ and $\e \sim \mathcal{N}(0,\sigma_e^{\star2})$, like in classical LMM defined in Equation \eqref{eq:LMM}. The heritability is then defined ``at the liability scale", which means for the continuous variable $\l$, and is given by the same expression \eqref{eq:eta_star_intro} as for quantitative traits.
\\
Several methods were established to estimate heritability in Model \eqref{eq:GLMM}: among them we can quote the MCMC method of \cite{hadfield} and the penalized quasi-likelihood approach of \citet{breslow}. The theoretical properties of these estimators have not been demonstrated and their numerical performances can be found in the comparative study of \citet{villemereuil2013}.
%
%
%
\\
Another modeling and definition for the heritability of  a binary trait, which is older and also more frequently used than the previous one, was proposed by \citet{falconer},  who assumed that the binary observations could be seen as an indicator function of a Gaussian variable exceeding a given threshold $t$:
\begin{equation}
\label{eq:liability_model}
\Y_i=\mathbbm{1}_{\{\l_i >t \}},
\end{equation}
\\
 with $\l_i$ defined by the same expression \eqref{eq:liability} than in Model \eqref{eq:GLMM}.
Observe that the threshold $t$ is directly linked to the prevalence of the disease in the population, that is the proportion $K$ of the population which is affected by the disease. Indeed,
\begin{equation}
 K =\mathbb{P}(\Y_i=1)= \mathbb{P}(\l_i>t).
 \label{eq:K_t}
\end{equation}
\\
%
The unobserved Gaussian variable $\l=(\l_1,\dots, \l_n)$ is also called the liability in this modeling, which is usually called the ``liability model" (\citet{falconer}, \citet{lee2011}, \citet{tenesa_2013}) and has been shown to be a reasonable modeling for complex diseases, for instance by \citet{purcell}. The heritability is then also defined at the liability scale as in Equation \eqref{eq:eta_star_intro}.
%
%
%
Regarding the procedures based on this modeling defined in Equations \eqref{eq:liability} and \eqref{eq:liability_model}, \citet{lee2011} proposed to use a maximum likelihood approach as if the binary traits were Gaussian, and then to apply a multiplicative factor to correct this approximation. \citet{golan2014}  showed that this heritability estimator was strongly biased in several realistic scenarios, in particular it was very sensitive to the prevalence of the disease (when the disease is rarer, the bias increases). The estimator also underestimates the heritability when the true heritability is high.
\citet{Weissbrod} introduced a maximum likelihood based strategy to rebuild the underlying liability before estimating the heritability. 
\\
However, all the aforementioned methods raise two main concerns: first, they have no theoretical validation. Second, they do not take into account an essential element of case-control studies: in a medical study, the number of patients is similar to the number of controls even though the studied disease might be rare, which means that the proportion of cases in the study does not reflect the proportion of cases in the population.
This oversampling of the cases has been noticed and handled by the approach of \citet{golan2014}, who proposed a moment based method to estimate the heritability.
The ground of their methodology was to compute an approximate quantity of the expectation $E$ of $\W_i \W_j$, for two individuals $i$ and $j$, $\W_i$ being a centered and normalized version of the binary data $\Y_i$, and conditionally to the fact that individuals $i$ and $j$ are in the study. This approach will be further described in Section \ref{sec:golan}.
\\
Since the method of \citet{golan2014},  presented very good numerical results but was not supported by theoretical grounds, we propose in this paper to investigate the theoretical properties of their method.
Our main result is to show that the least squares estimator obtained with the first order approximation of $E$ provides a consistent estimator of $\eta^\star$. We also propose a simulation study to compare the numerical performances of the estimators obtained with first and second order approximations of $E$. We show in particular that the computational times associated to the second order estimator are substantially larger with no obvious improvement from the statistical point of view.
\\
The model we study and the main definitions are given in Section \ref{sec:model}. Section \ref{sec:estimator} contains the first order approximation of the expectation $E$ with the corresponding estimator of $\eta^\star$ and Section \ref{sec:theorem} presents our consistency result for this estimator. The second order approximation of $E$ is given in Section \ref{sec:approx2} and the numerical comparison of the two estimators can be found in Section \ref{sec:numerical}. In Section \ref{sec:discussion}, we discuss the results and potential perspectives. Finally, the proofs are given in Section \ref{sec:proofs}.

\section{Model and definitions}
\label{sec:model}
\subsection{Liability model}
Let us denote $K$ the prevalence of a disease in a population, that is the proportion of the population affected by the disease.
Let $\Y_i$ be the random variable such that $\Y_i=1$ if the individual $i$ is ill (then, individual $i$ is called a case) and $\Y_i=0$ if the individual $i$ is healthy (then individual $i$ is called a control). We assume that the $\Y_i$'s are linked to unobserved variables $\l_i$ as follows 
\begin{equation}
\Y_i=\mathbbm{1}_{\{\l_i >t \}},
\label{eq:liability_model2}
\end{equation}
\\
 where $t$ is a given threshold, related to the prevalence $K$ by \eqref{eq:K_t}, and the $\l_i$'s are defined as

\begin{equation}
\l=\Z \u + \e,
\label{eq:liability2}
\end{equation}
\\
where $\l=(\l_1,\dots, \l_n)$, $\u$ and $\e$ are random effects such that $\u \sim \mathcal{N}(0,\sigma_u^{\star 2} \Id_{\rset^N})$ and 
$\e \sim \mathcal{N}(0,\sigma_e^{\star 2} \Id_{\rset^n})$. 
The vector $\u$ corresponds to the genetic effects and $\e$ to the environmental effects. Moreover, $\Z$ is a $n\times N$ random
matrix which contains the genetic information, and which is such that the $\Z_{i,k}$ are normalized random variables in the following sense: they are defined from 
a matrix $A=(A_{i,k})_{1\leq i\leq n,\,1\leq k\leq N}$ by
\begin{equation}\label{eq:normalization_1}
\Z_{i,k}=\frac{A_{i,k}-\overline{A}_{k}}{s_{k}},\;i=1,\dots,n,\;k=1,\dots,N\;,
\end{equation}
where
\begin{equation}\label{eq:normalization_2}
\overline{A}_{k}=\frac{1}{n}\sum_{i=1}^{n}A_{i,k},\;s_{k}^{2}=\frac{1}{n}\sum_{i=1}^{n}(A_{i,k}-\overline{A}_{k})^{2},
\;k=1,\ldots,N\;.
\end{equation}
In (\ref{eq:normalization_1}) and (\ref{eq:normalization_2}) the $A_{i,k}$'s are such that
for each $k$ in $\{1,\dots,N\}$ the $(A_{i,k})_{1\leq i\leq n}$ are independent and identically distributed 
random variables and such that the columns of $A$ are independent. 
\\
In practice, the matrix $A$ contains the genetic information about all the individuals in the study. More precisely, for each $k$, $A_{i,k}=0$ (resp. 1, resp. 2) if the genotype of the $i$th individual at locus $k$ is qq
(resp. Qq, resp. QQ). In this paper, we consider a more general case with mild assumptions on the distribution of the random variables $A_{i,k}$, which are described in Section \ref{sec:theorem}.
\\
With the definition \eqref{eq:normalization_1}, the columns of $\Z$ are empirically centered and normalized, and 
one can observe that
$$
\Var(\l|\Z)=N{\sigma_u^\star}^2 \R +{\sigma_e^\star}^2 \Id_{\rset^n}\;,\textrm{ where }\R=\frac{\Z\Z'}{N}.
$$
\\
The heritability at the liability scale, which is the parameter we want to estimate, is defined as the ratio of variances: 
\begin{equation}
\eta^\star=\frac{N\sigma_u^{\star2}}{N\sigma_u^{\star2}+\sigma_e^{\star2}}.
\label{eq:eta_star}
\end{equation}
\\
The variance of $\l$ conditionally to $\Z$ can then be rewritten with respect to $\eta^\star$ and $\sigma^{\star2}=N\sigma_u^{\star2}+\sigma_e^{\star2}$ as: 
$$
\Var(\l|\Z)=\eta^\star{\sigma^\star}^2 \R+(1-\eta^\star){\sigma^\star}^2 \Id_{\rset^n}\;.
$$
We will assume in the sequel without loss of generality that $\sigma^{\star2}=1$. Indeed, if $\sigma^{\star2} \neq 1$, we can consider the variable $\l'_i=\frac{\l_i}{\sigma^{\star}}$ and 
then, instead of estimating $t$ from the prevalence $K$ with the relationship \eqref{eq:K_t}, we estimate directly $t/\sigma^{\star}$.

\subsection{Case control study}
Since the prevalence $P$ in the study can be very different from the prevalence $K$ in the general population (the cases are substantially oversampled in a case-control study), 
it is essential to consider that the observations that we have access to depend on the probabilities for both cases and controls to be selected in the study.
Indeed, if $p_{control}$ denotes the probability for a control to be selected in the study, we can define the corresponding variable $U_i\sim \mathcal{B}(0,p_{control})$ which is equal to $1$ if individual $i$ is part of the study. Similarly we define the probability $p_{case}$ for a case to be selected for the study and the corresponding variable $V_i\sim \mathcal{B}(0,p_{case})$.
Then for any individual $i$, we define the variable $\epsilon_i$ by 
 $$\epsilon_i= V_i \Y_i + U_i (1-\Y_i),$$
 which is equal to $1$ if individual $i$ belongs to the study and $0$ if not. We assume that the variables $U_1, \dots, U_n, V_1, \dots, V_n$ are independent and independent of $\Y_1, \dots, \Y_n$ and $\Z$.
\\
Since we do not observe $\Y_i$ for the whole population but only for the individuals who belong to the study, we will work with the variables $\W_i$ defined by
$$\W_i= \frac{\Y_i-P}{\sqrt{P(1-P)}} \epsilon_i,$$
\\
which are centered versions of $\Y_i$ in the study and are non-zero only if individual $i$ belongs to the study.
\\
The probabilities $p_{case}$ and $p_{control}$ are chosen such that the prevalence in the study is equal to $P$. Indeed, if we assume that
\begin{equation}
p_{case}=1,
\label{eq:full_asc}
\end{equation}
it implies that 
 \begin{equation}
 p_{control}= \frac{K(1-P)}{P(1-K)}.
 \label{eq:p_control}
 \end{equation}
 \\
The proof of \eqref{eq:p_control} is given in Appendix \ref{app3}.
Equation \eqref{eq:full_asc} means that all cases are accepted in the study and it is usually called a ``full ascertainment" assumption (see for instance \cite{golan2014}).
\\
%
\section{Heritability estimator}
\label{sec:estimator}

\subsection{Method of \cite{golan2014}}
\label{sec:golan}

\cite{golan2014} considered a simplified version of Model \eqref{eq:liability}, where the liability is given by

$$\l= \g + \e,$$
\\
where $\g$ is a genetic random effect, which can be correlated across individuals, and $\e$ is the environmental
random effect, which is assumed to be independent of the genetic effect. Both effects
are assumed to be Gaussian: $\e$ has a variance equal to $(1 -\eta^\star) \Id_{\rset^n} $ and $\g$ has a covariance matrix, the diagonal entries of which are equal to $\eta^\star$ and the  non diagonal term $(i,j)$ is equal to $\eta^{\star} \G_{i,j}$.  The covariance matrix of $(\l_i,\l_j)$ is given by
$$\Sigma=\begin{pmatrix}
   1 & \eta^\star\G_{i,j} \\
   \eta^\star\G_{i,j}  & 1
\end{pmatrix}.$$
\\
The heritability estimator proposed by \cite{golan2014} is a least square estimator obtained by minimizing 

\begin{equation}
\underset{i \neq j}{\sum} \left( \W_i \W_j - \mathbb{E}[\W_i \W_j \vert  \epsilon_i=\epsilon_j=1] \right)^2.
\label{eq:LS}
\end{equation}
\\
Since the expression of $\mathbb{E}[\W_i \W_j \vert \epsilon_i=\epsilon_j=1]$ has no explicit formula as we shall see hereafter, \cite{golan2014} proposed to take advantage of the fact that the correlations $\G_{i,j}$ are small for $i \neq j$. 
\\
The ground of the method is to write 
\\
\begin{align}
\label{E_Wi_Wj}
\mathbb{E}[\W_i \W_j \vert \epsilon_i=\epsilon_j=1] = \frac{  \frac{1-P}{P}  \mathbb{P}(\Y_i=\Y_j=1 )   - \frac{K(1-P)}{P(1-K)}\mathbb{P}(\Y_i \neq \Y_j ) + \frac{K^2(1-P)}{P(1-K)^2}\mathbb{P}(\Y_i=\Y_j=0 ) }{\mathbb{P}(\Y_i=\Y_j=1 ) + \left(\frac{K(1-P)}{P(1-K)}\right)^2\mathbb{P}(\Y_i=\Y_j=0 ) + \frac{K(1-P)}{P(1-K)}\mathbb{P}(\Y_i \neq \Y_j) }
\end{align}
\\
and to propose approximations of $\mathbb{P}(\Y_i \neq \Y_j)$, $\mathbb{P}(\Y_i=\Y_j=0)$ and $\mathbb{P}(\Y_i \neq \Y_j)$ thanks to Taylor developments around the quantity $\G_{i,j}$.
The computations leading to \eqref{E_Wi_Wj} can be found in Appendix \ref{app1}.
\\
This approximation, plugged in the least squares criterion \eqref{eq:LS}, led to the heritability estimator given by 
\begin{align}
\hat{\eta}=\left[\frac{\underset{i \neq j}{\sum} \W_i \W_j \G_{i,j} }{c \underset{i \neq j}{\sum}  \G_{i,j}^2} \wedge 1 \right] \vee 0,
\label{eq:estim_golan}
\end{align}
\\
 where 
\\
\begin{align}
 c= \phi(t)^2 \frac{P(1-P)}{K^2(1-K)^2},
 \label{eq:c}
 \end{align} 
 \\
  $\phi$ being the density of the standard Gaussian distribution.
%
%

%
%

\subsection{Our method}
In Model \eqref{eq:liability_model2} that
we consider, 
the variance matrix $\Sigma^{(N)}$ of $(\l_i,\l_j)$ conditionally to $\Z$ can be written as 

\begin{align*}
\Sigma^{(N)}=  \begin{pmatrix}
   1 +  \eta^\star (\G_N(i,i)-1)& \eta^\star\G_N(i,j) \\
  \eta^\star \G_N(i,j) & 1 + \eta^\star (\G_N(j,j)-1)
\end{pmatrix},
\end{align*} 
where for all $ 1 \leq i,j \leq n$, 

\begin{equation}
\G_N(i,j)= \frac{1}{N} \overset{N}{\underset{k=1}{\sum}} \Z_{i,k} \Z_{j,k}.
\label{eq:G_ij}
\end{equation}
Note that in the model we consider, $\G_N(i,j)$ is a random variable, which is not the case of the quantity $\G_{i,j}$ in the model studied by \cite{golan2014}.
A key element is to notice that $\Sigma^{(N)}$ is close to the $n \times n$ identity matrix, more precisely 
\begin{align}
\Sigma^{(N)}=  \begin{pmatrix}
   1 +  \eta^\star\frac{A_N(i)}{\sqrt{N}}& \eta^\star\frac{B_N(i,j)}{\sqrt{N}} \\
  \eta^\star \frac{B_N(i,j)}{\sqrt{N}} & 1 + \eta^\star \frac{A_N(j)}{\sqrt{N}} 
\end{pmatrix} 
\label{eq:Sigma_N}
\end{align} 
\\
where $A_N(i)=O_p(1)$, $A_N(j)=O_p(1)$ and $B_N(i,j)=O_p(1)$. The proof of \eqref{eq:Sigma_N} can be found in Appendix \ref{app2}.
\\
Then, following the idea of \cite{golan2014}, we propose to approximate $$\mathbb{E}[\W_i \W_j \vert \Z, \epsilon_i=\epsilon_j=1]$$ defined in Equation \eqref{E_Wi_Wj} thanks to Taylor developments around $\frac{A_N(i)}{\sqrt{N}}$, $\frac{A_N(j)}{\sqrt{N}}$ and $\frac{B_N(i,j)}{\sqrt{N}}$.
\\
The detailed computations are devised in Section \ref{sec:taylor}.
\\
A first order approximation of $\mathbb{E}[\W_i \W_j \vert \Z, \epsilon_i=\epsilon_j=1]$, plugged in \eqref{eq:LS}, leads to the same estimator $\hat{\eta}^{(1)}$ as the one proposed by \cite{golan2014}. Indeed, we obtain
\\
\begin{align}
\hat{\eta}^{(1)}=\left[\frac{\underset{i \neq j}{\sum} \W_i \W_j \G_N(i,j) }{c \underset{i \neq j}{\sum}  \G_N(i,j)^2} \wedge 1 \right] \vee 0,
\label{eq:estim_nous}
\end{align}
\\
 where 
 $c= \phi(t)^2 \frac{P(1-P)}{K^2(1-K)^2}.$
 \\
 In Section \ref{sec:approx2}, we consider the second order approximation, which is different from the one devised by \cite{golan2014}.

$$ $$

 \section{Consistency of the heritability estimator $\hat{\eta}^{(1)}$}
\label{sec:theorem}

In this section, we consider the heritability estimator $\hat{\eta}^{(1)}$ defined in Equation \eqref{eq:estim_nous}.
\\
\begin{assumption}
There exist $d>0$, $C>0$ and a neighborhood $V_0$ of $0$ such that for all $\lambda$ in $V_0$

\begin{enumerate}[label=\textbf{1.\arabic*}]
\item \label{hyp1} $\mathbb{E}[\exp \left( \lambda (A_{i,k}-\mathbb{E}[A_{i,k}])^2-\sigma_k^2\right)] \leq C \exp (d\lambda^2)$ 
\item \label{hyp2} $\mathbb{E}[\exp \left(\lambda  ( A_{i,k}-\mathbb{E}[A_{i,k}])\right)] \leq C \exp (d\lambda^2)$ 
\item \label{hyp4}$\mathbb{E}[\exp \left(\lambda  ( A_{i,k}-\mathbb{E}[A_{i,k}])( A_{j,k}-\mathbb{E}[A_{i,k}])\right)] \leq C \exp (d\lambda^2)$
\end{enumerate}
for all $i \neq j$ and for all $k$, where the $A_{i,k}$'s are defined in \eqref{eq:normalization_1} and $\sigma_k^2$ is the variance of $A_{i,k}$.
\label{HYP:HYP1}
\end{assumption}

\begin{assumption}
$$ $$
\vspace{-10mm}
\begin{enumerate}[label=\textbf{2.\arabic*}]
\item \label{hyp3}
 $\underset{k=1..N}{\inf} \sigma_k^2=\delta_{min} >0$
\item \label{hyp5} $\underset{k=1..N}{\sup} \sigma_k^2=\delta_{max} < +\infty$
\end{enumerate}
\label{HYP:HYP2}
\end{assumption}

\begin{theorem}

Let $\Y=(\Y_1,…,\Y_n)$ satisfy Model \eqref{eq:liability_model} with $A$ satisfying Assumptions \ref{HYP:HYP1} and \ref{HYP:HYP2}, and $\hat{\eta}^{(1)}$ the estimator of $\eta^\star$ defined in Equation \eqref{eq:estim_nous}. Then, as $n,N \to \infty$ such that $n/N \to a \in (0,+\infty)$, 
$$ \hat{\eta}^{(1)}=\eta^\star + o_p(1).$$
\label{th:consistency}
\end{theorem}
The proof of Theorem \ref{th:consistency} relies on the following lemmas.

\begin{lemma}
When $n$ and $N$ go to infinity and $n/N$ goes to $a$,
 $$\frac{1}{n}\underset{i \neq j}{\sum}  \G_N(i,j)^2 \textrm{converges in probability to $a$.}$$ 
 \label{lemma:denominator}
\end{lemma}
We will then have to focus on 
\\
\begin{align}
\frac{1}{n}\underset{i \neq j}{\sum} \W_i \W_j \G_N(i,j) & =  \left[\frac{1}{n}\underset{i \neq j}{\sum}( \W_i \W_j - \mathbb{E}[\W_i \W_j \vert \Z, \epsilon_i=\epsilon_j=1])\G_N(i,j) \right. \nonumber \\
& +\left. \frac{1}{n}\underset{i \neq j}{\sum} \mathbb{E}[\W_i \W_j \vert \Z, \epsilon_i=\epsilon_j=1]\G_N(i,j) \right].
 \label{eq:decomposition}
\end{align}
\\
Let $E_N$ be the following event

$$E_N=\left\{\underset{i}{ \sup} \vert \G_N(i,i)-1 \vert \leq \epsilon_N \textrm{ and } \underset{i \neq j}{ \sup} \vert \G_N(i,j) \vert \leq \epsilon_N\right\},$$
\\
where $\epsilon_N=\frac{1}{N^{\frac{1}{2}-\gamma}}$ with $\gamma$ a positive number such that $\gamma < 1/10$.
\\
Let us denote $E_N^c$ the complement of the event $E_N$. We consider the following decomposition  

$$\hat{\eta}^{(1)}=\hat{\eta}^{(1)} \mathbbm{1}_{E_N} + \hat{\eta}^{(1)} \mathbbm{1}_{E_N^c }.$$

\begin{lemma}
 For all values of $q$, the probability of $E_N^c$ satisfies 
$\mathbb{P}(E_N^c)=O(\frac{1}{N^q})$ when $ N \to +\infty$.
\label{lemma5}
\end{lemma}

Using the result of Lemma \ref{lemma5}, $\hat{\eta}^{(1)} \mathbbm{1}_{E_N^c }$ converges in probability to $0$ since 

$$\mathbb{E}[\vert \hat{\eta}^{(1)} \mathbbm{1}_{E_N^c } \vert ] \leq \mathbb{P}(E_N^c)=O\left(\frac{1}{N^q}\right).$$

\begin{lemma}
When $n$ and $N$ go to infinity and $n/N$ goes to $a\in (0,+\infty)$,
$$\frac{1}{n}\underset{i \neq j}{\sum} \mathbb{E}[\W_i \W_j \vert \Z, \epsilon_i=\epsilon_j=1]\G_N(i,j) \mathbbm{1}_{E_N}$$

converges in probability to $a c \eta^\star$, where $c$ is defined in Equation \eqref{eq:c}.
\label{lemma3}
\end{lemma}

\begin{lemma}
When $n$ and $N$ go to infinity and $n/N$ goes to $a\in (0,+\infty)$,
$$\frac{1}{n}\underset{i \neq j}{\sum}( \W_i \W_j - \mathbb{E}[\W_i \W_j \vert \Z,\epsilon_i=\epsilon_j=1])\G_N(i,j)\mathbbm{1}_{E_N}$$

converges in probability to $0$.
\label{lemma4}
\end{lemma}

The results of Lemmas \ref{lemma3} and \ref{lemma4} achieve the proof of Theorem \ref{th:consistency}. 
\\
The proof of Lemmas \ref{lemma:denominator}, \ref{lemma5}, \ref{lemma3} and \ref{lemma4} are given in Section \ref{sec:proof_consistency}.

 \section{Second order approximation of  $\mathbb{E}[\W_i \W_j \vert \Z, \epsilon_i=\epsilon_j=1]$}
 \label{sec:approx2}
 
The purpose of this section is to study the behaviour of the heritability estimator obtained thanks to a second order approximation of $\mathbb{E}[\W_i \W_j \vert \Z, \epsilon_i=\epsilon_j=1]$.
\\
Instead of computing the approximation till order $1/\sqrt N$, we compute the approximation till order $1/N$ and we obtain:
\\
      \begin{align*}
\mathbb{E}[\W_i \W_j \vert & \Z, \epsilon_i=\epsilon_j=1]= \frac{\eta^\star}{\sqrt N} \frac{P(1-P)}{K^2(1-K)^2}\phi(t)^2B_N(i,j) +\frac{t^2}{4} \frac{\eta^{\star 2}}{N}A_N(i) A_N(j) \frac{P(1-P)}{K^2(1-K)^2} \\
&+  \frac{\eta^{\star 2}}{N} \frac{P(1-P)}{K^2(1-K)^2}\phi(t)^2B_N(i,j)^2 \left[\frac{t^2}{2}-\frac{(P-K)^2}{K^2(1-K)^2}\right] \\
&+ \frac{\eta^{\star 2}}{2 N} \frac{P(1-P)}{K^2(1-K)^2}\phi(t)^2B_N(i,j)(A_N(i)+A_N(j)) \left[t^2-1 - \frac{P-K}{K(1-K)}t\phi(t)\right]+ O_p\left(\frac{1}{N^{\frac{3}{2}}}\right)
  \end{align*}
\\ 
  The proof of this computation is detailed in Section \ref{pr:approx2}.
\\ 
 Since the minimizer in $\eta$ of the quantity
\\  
\begin{align*}
g(\eta)= &  \underset{i \neq j}{\sum} \left( \W_i \W_j - \frac{\eta}{\sqrt N} \frac{P(1-P)}{K^2(1-K)^2}\phi(t)^2B_N(i,j) -\frac{t^2}{4} \frac{\eta^{ 2}}{N}A_N(i) A_N(j) \frac{P(1-P)}{K^2(1-K)^2} \right. \\
&- \frac{\eta^{2}}{N} \frac{P(1-P)}{K^2(1-K)^2}\phi(t)^2B_N(i,j)^2 \left[\frac{t^2}{2}-\frac{(P-K)^2}{K^2(1-K)^2}\right] \\
& \left.- \frac{\eta^{ 2}}{2 N} \frac{P(1-P)}{K^2(1-K)^2}\phi(t)^2B_N(i,j)(A_N(i)+A_N(j)) \left[t^2-1 - \frac{P-K}{K(1-K)}t\phi(t)\right]\right)^2
\end{align*}
 \\
 has no explicit form, we use a Newton-Raphson approach to obtain the corresponding heritability estimator $\hat{\eta}^{(2)}$ of the second order approximation.
\\Note that the second order approximation, which depends on $B_N(i,j)$ but also on $A_N(i)$ and $A_N(j)$, is different from the one found by \citet{golan2014}. 

 \section{Numerical study}
\label{sec:numerical}

In this section, we propose to study the numerical performance of the estimators $\hat{\eta}^{(1)}$ and $\hat{\eta}^{(2)}$ devised respectively in Sections \ref{sec:estimator} and \ref{sec:approx2}. Since \citet{golan2014} already compared the estimator $\hat{\eta}^{(1)}$ to the one proposed by \citet{lee2011} and stated several arguments in favor of their estimator, we will focus on comparing our two estimators in terms of statistical and computational efficiency. 

\subsection{Simulation process}

 In this simulation study, we generated data sets with $n \simeq 200$, $N=10000$ in order to respect the classical scenario where $N >>n$. The value of the prevalence in the population varies from $0.005$ to $0.1$. The observations were generated as follows. 
 
\begin{itemize}
\item We set the parameters $\eta^\star$, $K$, $P=1/2$ and the size of the general population, chosen very large.  Notice that the number of individuals selected in the study varies from one sample to another. We chosed in practice a population size in order to have around $100$ patients in the study.

\item We generated the Gaussian random effects $\u$ and $\e$ with respective variances
$\sigma_u^{\star 2} = \eta^\star/N$ and $\sigma_e^{\star2}=1-\eta^\star$.

\item We generated liabilities, from which we generated binary observations in order to have a prevalence equal to $K$ in the general population. 

\item For each individual, we determined those who stayed in the study: the cases are automatically selected (full ascertainment assumption) but each control is selected with probability $p_{control}$ computed in Equation \eqref{eq:p_control}. 
\end{itemize}

\subsection{Results}
 \begin{figure}[!h]

\begin{center}

\begin{tabular}{cc}

$\eta^\star=0.5$ & $\eta^\star=0.7$ \\
\includegraphics[width=0.35\textwidth,trim=10mm 15mm 10mm 20mm,clip=true]{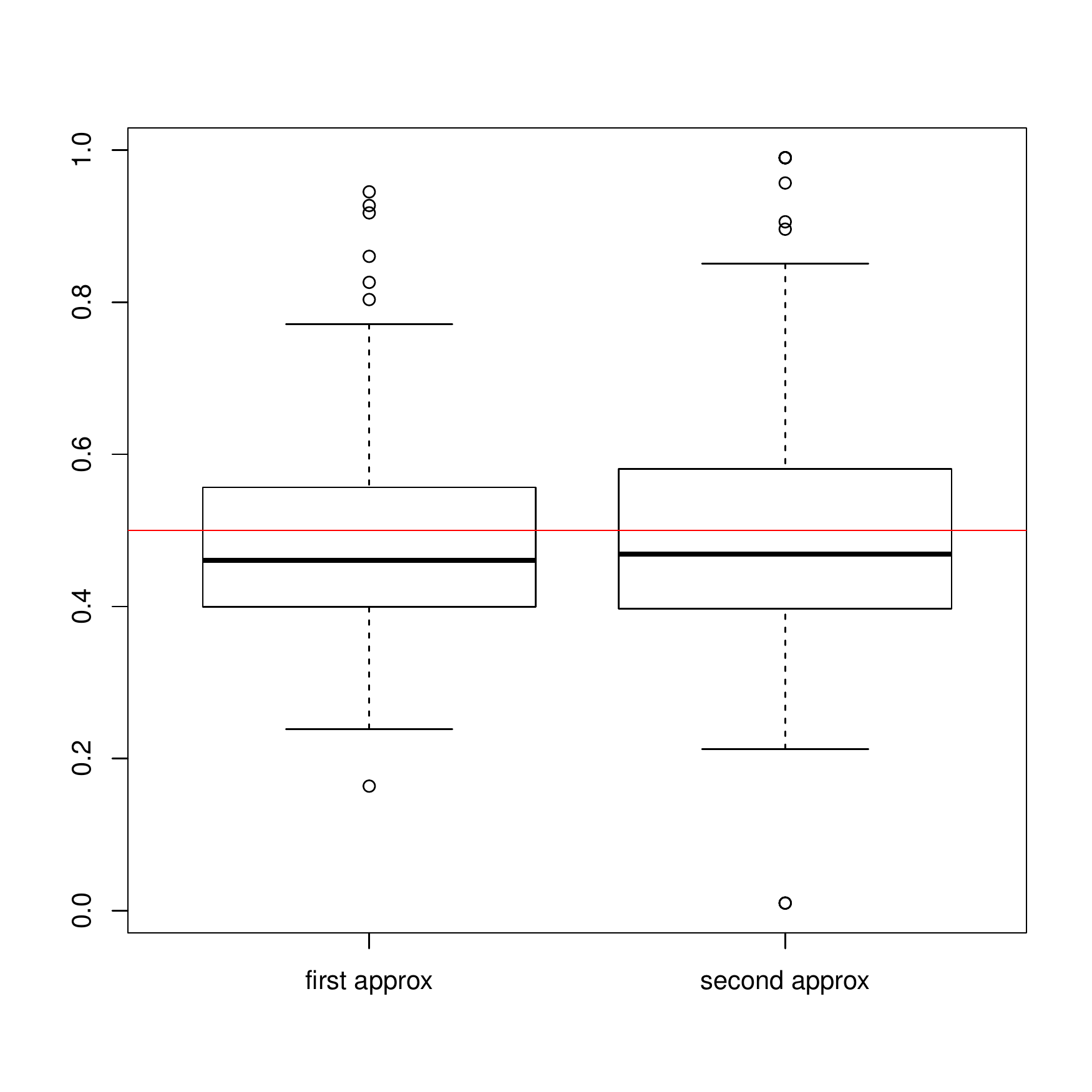} &

\includegraphics[width=0.35\textwidth,trim=10mm 15mm 10mm 20mm,clip=true]{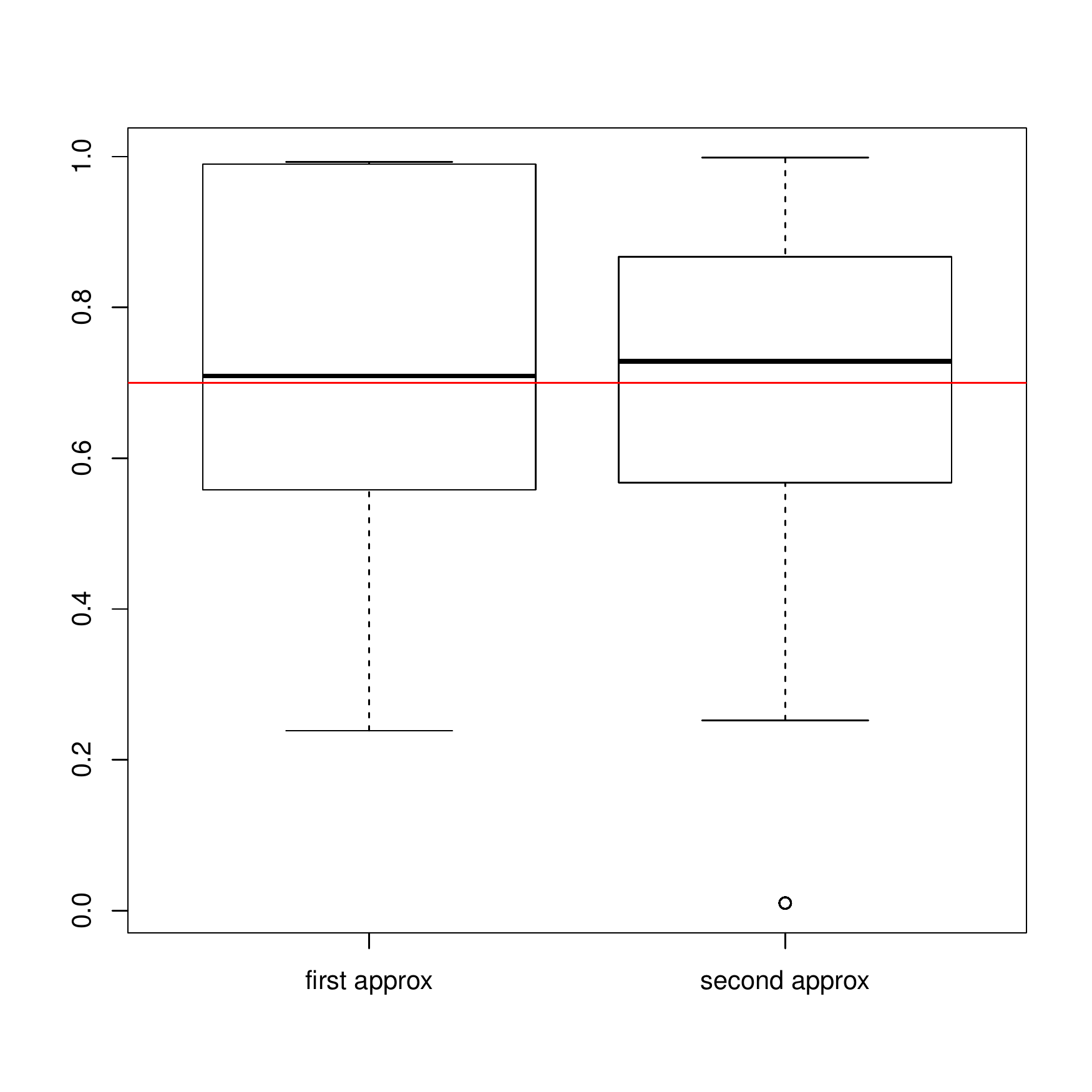} \\

\includegraphics[width=0.35\textwidth,trim=10mm 15mm 10mm 20mm,clip=true]{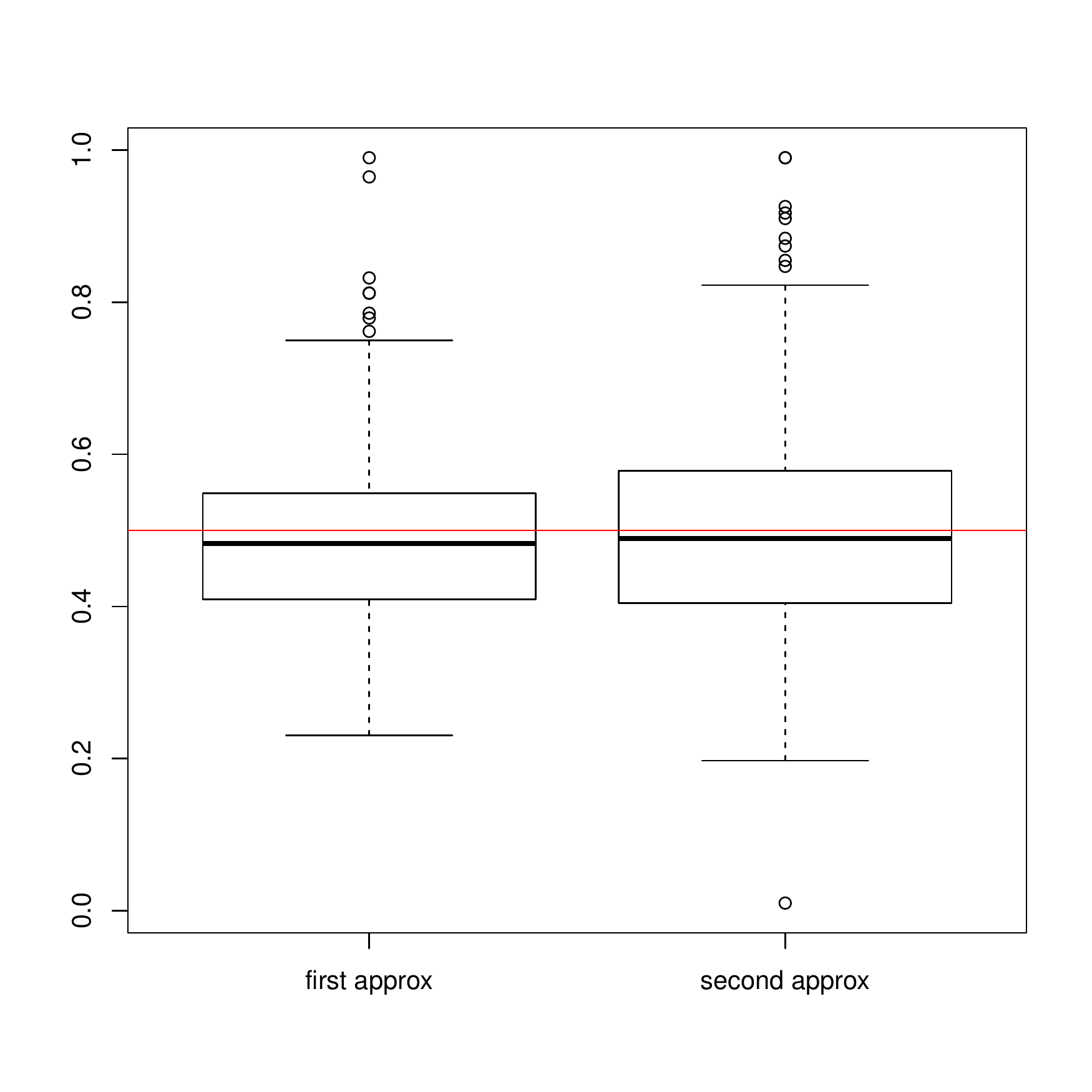} &

\includegraphics[width=0.35\textwidth,trim=10mm 15mm 10mm 20mm,clip=true]{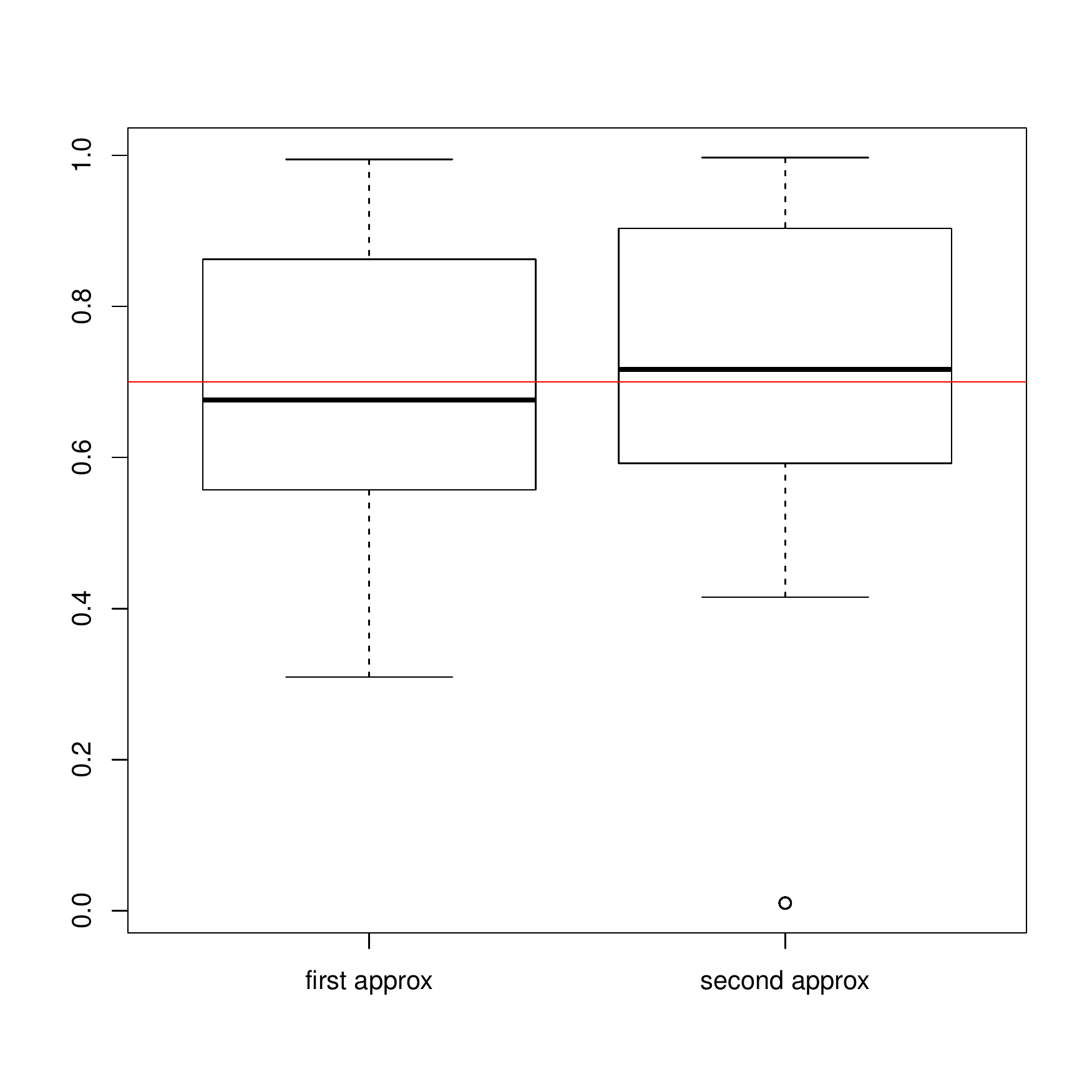} \\

\includegraphics[width=0.35\textwidth,trim=10mm 15mm 10mm 20mm,clip=true]{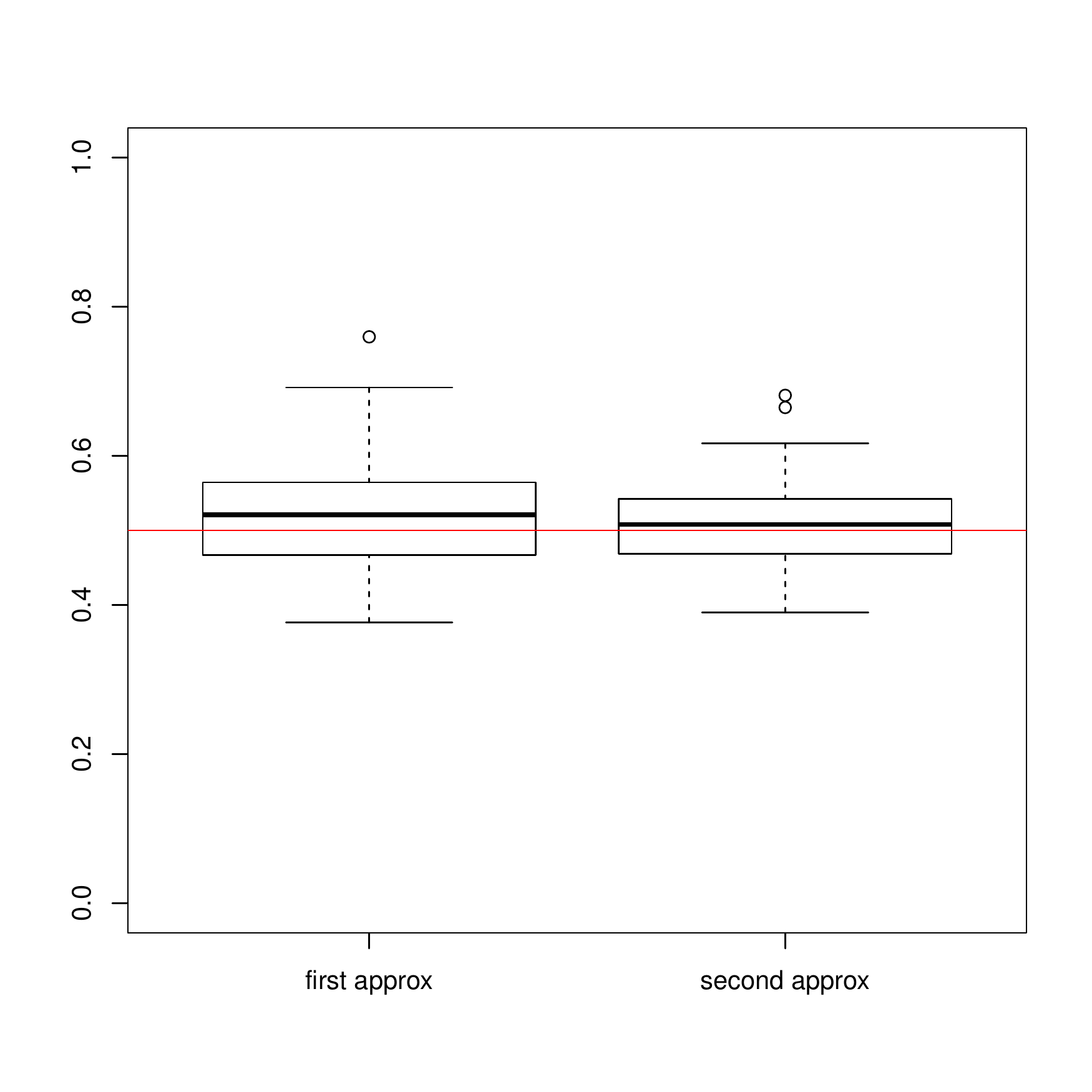} &

\includegraphics[width=0.35\textwidth,trim=10mm 15mm 10mm 20mm,clip=true]{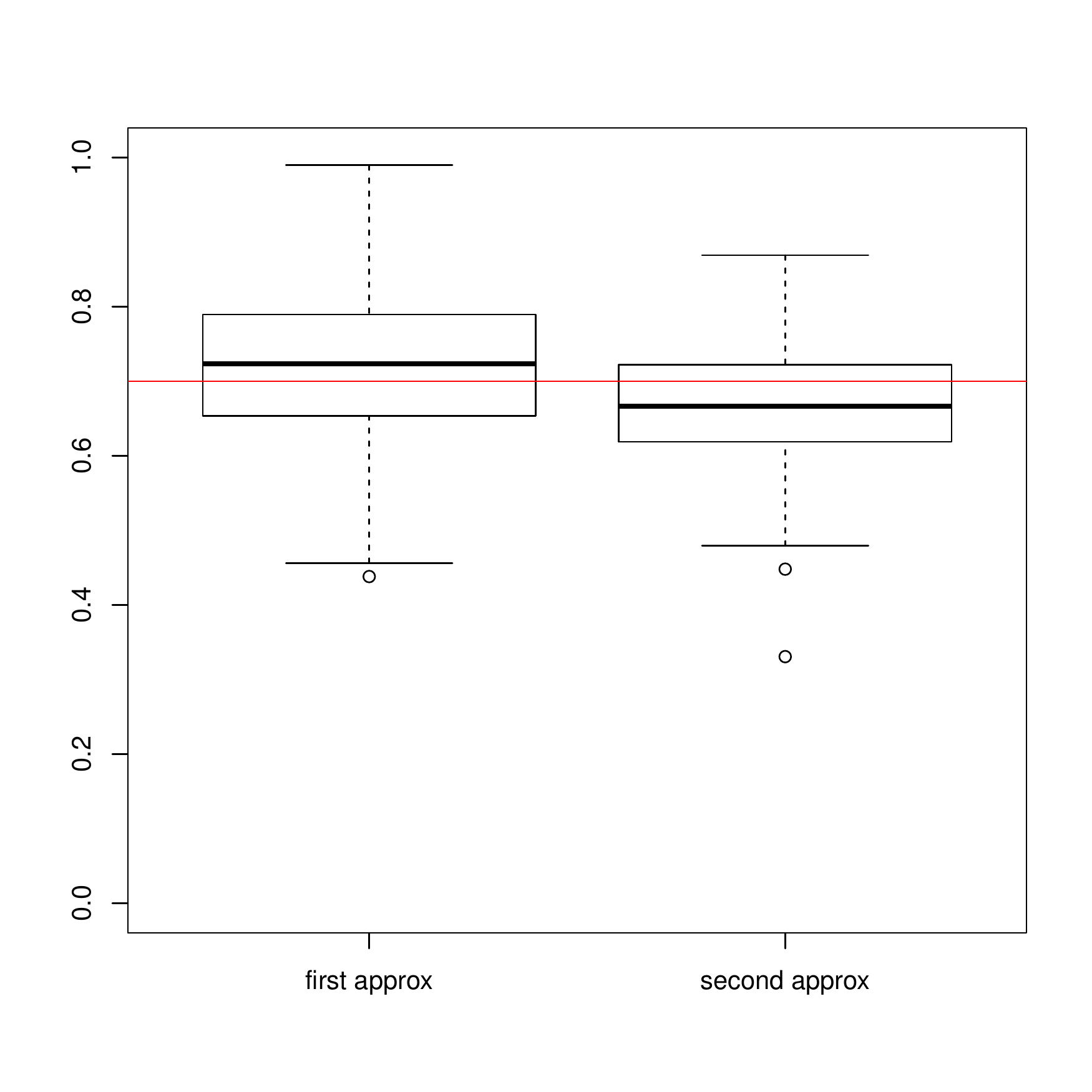}
\end{tabular}

\caption{Boxplots for $\hat{\eta}^{(1)}$ (``first approx") and $\hat{\eta}^{(2)}$ (``second approx") for different values of $\eta^\star$: $0.5$ (left), $0.7$ (right) and different values of the prevalence $K$: $0.005$ (top), $0.01$ (middle), $0.1$ (bottom). The sample size is $n \simeq 200$ and $N=10000$. Each boxplot is generated from $200$ replications.}
\label{fig:box_kvar}
\end{center}
\end{figure} 

\begin{figure}[!h]

\begin{center}

\begin{tabular}{cc}

$n \simeq 100$ & $n \simeq 1000$ \\
\includegraphics[width=0.35\textwidth,trim=10mm 15mm 10mm 20mm,clip=true]{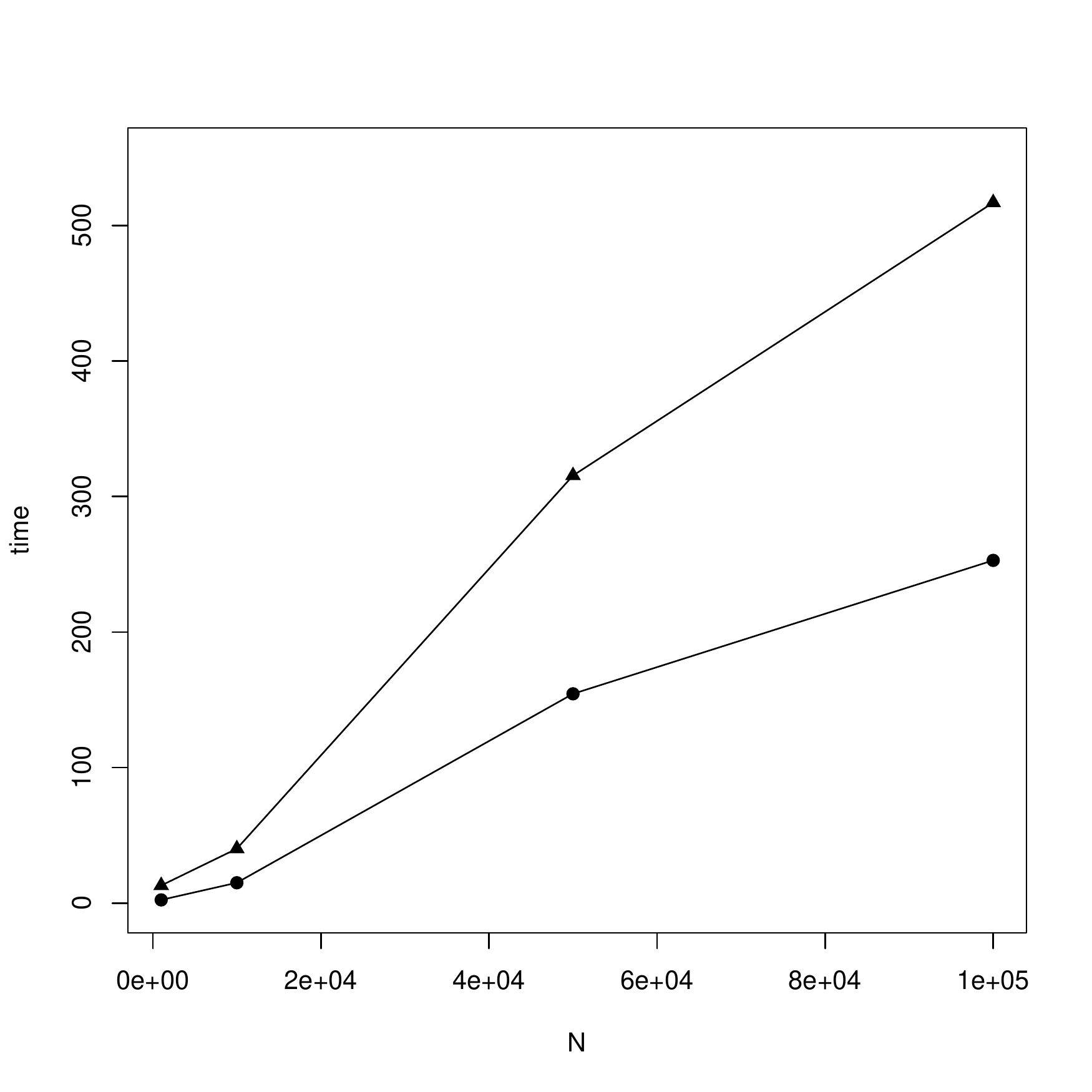} &

\includegraphics[width=0.35\textwidth,trim=10mm 15mm 10mm 20mm,clip=true]{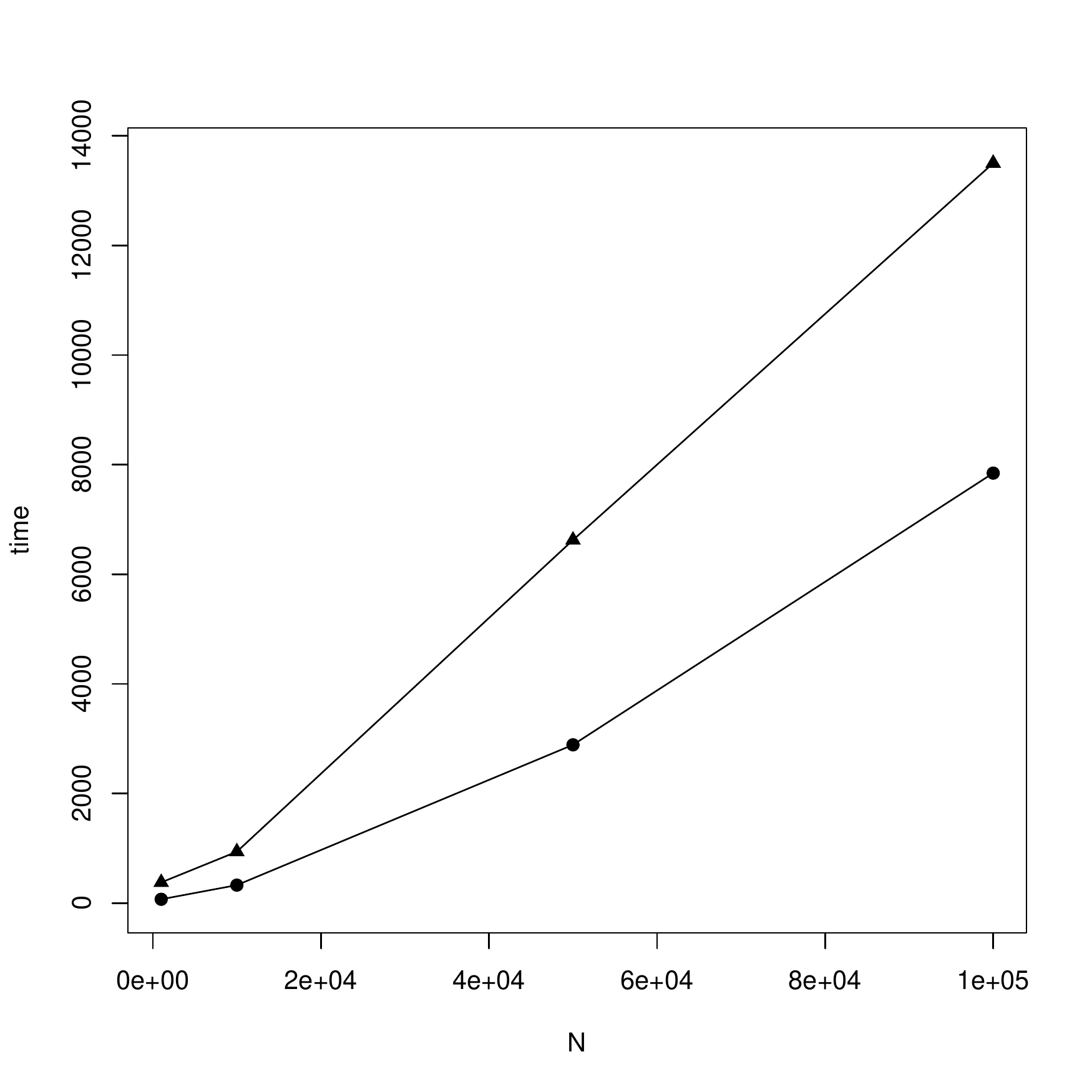} 
\end{tabular}

\caption{Time in seconds to compute an estimation of $\eta^\star$ obtained with $\hat{\eta}^{(1)}$ and $\hat{\eta}^{(2)}$ for $n \simeq 100$ (left) and $n \simeq 1000$ (right) and for different values of $N$ (from $1000$ to $10^5$).}
\label{fig:temps}
\end{center}
\end{figure}

Figure \ref{fig:box_kvar} displays the estimations of $\eta^\star$ obtained with both estimators $\hat{\eta}^{(1)}$ and $\hat{\eta}^{(2)}$. First, we can notice that both estimators seem empirically unbiased. Second, we observe no obvious improvement of the performance of $\hat{\eta}^{(2)}$ compared to $\hat{\eta}^{(1)}$ in terms of empirical variance.  Finally, we can also note that the estimations seem more accurate when the prevalence $K$ is high, namely $K=0.1$.
\\
Table \ref{tab:comp_time} and Figure \ref{fig:temps} show the computational performance of both estimators. The computation of the estimator $\hat{\eta}^{(2)}$ obtained with the more refined approximation is obviously slower, but for small values of $n$ (namely, $n=100$), the time required to compute an estimation of $\eta^\star$ remains quite small (86 seconds, against 40 seconds for the other estimator). However, when $n$ is larger, the computational time increases substantially and the ``slower" estimator needs up to 13500 seconds, that is almost 4 hours, to compute an estimation of $\eta^\star$. 
\\
In conclusion, both estimators are empirically unbiased and since the computation of the estimator $\hat{\eta}^{(2)}$ is slower and does not improve the estimations of $\eta^\star$, we are satisfied with the first order approximation and the corresponding estimator $\hat{\eta}^{(1)}$.


 \begin{table}[!h]
 \centering
   \caption{Times in seconds to compute an estimation of $\eta^\star$ obtained with $\hat{\eta}^{(1)}$ and $\hat{\eta}^{(2)}$ for different values of $n$ ($100$ and $1000$) and $N$ (from $1000$ to $10^5$).}
\begin{tabular}{| c | c | c | c | c | c |} 
  \hline
 n & $N$ & $1000$ & $10000$ & $50000$ & $10^5$  \\
  \hline
 100 & $\hat{\eta}^{(1)}$ & 0.478  & 2.390 & 28.595 & 40.528  \\
  \cline{2-6}   

   & $\hat{\eta}^{(2)}$ & 3.148  & 7.127 & 56.761  & 86.156  \\
  \hline
  1000 & $\hat{\eta}^{(1)}$ & 69.047 & 327.240  & 2887.518  & 7845.281 \\

  \cline{2-6}  
   & $\hat{\eta}^{(2)}$ & 376.363 & 936.845 & 6624.186   &13500.510  \\
  \hline
  
\end{tabular}
\label{tab:comp_time}
\end{table}

\section{Discussion}
\label{sec:discussion}
In this paper, we proposed theoretical grounds to support the heritability estimator in case-control studies developed by \cite{golan2014}. We proved indeed its consistency in the framework where both the number of individuals $n$ and the number $N$ of SNPs go to infinity, when the ratio $n/N$ goes to a constant $a$.
\\
It would be interesting to complete this work with theoretical results which could allow the user to compute accurate confidence intervals, similarly to existing results for quantitive traits. As it is often the case in genetic applications, the question of removing strong assumptions such as the Gaussianity of the random effects or the independence of the columns of the SNP matrix remains a challenging issue. Considering possible sparsity in the random effects would also be an interesting improvement and will be the subject of a future work.

\section{Proofs}
\label{sec:proofs}

\subsection{Taylor development of $\mathbb{E}(\W_i \W_j \vert \Z, \epsilon_i=\epsilon_j=1)$ in Model \eqref{eq:liability_model2}}
\label{sec:taylor}

According to Equation \eqref{E_Wi_Wj}, we only need to compute approximations of $\mathbb{P}(\Y_i=\Y_j=1 \vert \Z)$, $\mathbb{P}(\Y_i=\Y_j=0 \vert \Z)$ and $\mathbb{P}(\Y_i \neq \Y_j \vert \Z)$ to obtain an approximation of  $\mathbb{E}[\W_i \W_j \vert \Z, \epsilon_i=\epsilon_j=1]$.

\begin{align*}
\mathbb{P}(\Y_i=\Y_j=1 \vert \Z)= \int_{t}^{\infty} \int_{t}^{\infty} f(x,y) dx dy,
\end{align*}

\begin{align*}
\mathbb{P}(\Y_i=\Y_j=0 \vert \Z)= \int_{-\infty}^{t} \int_{-\infty}^{t} f(x,y) dx dy
\end{align*}
\\
 and \begin{align*}
\mathbb{P}(\Y_i \neq \Y_j \vert \Z)= 2 \int_{-\infty}^{t} \int_{t}^{\infty} f(x,y) dx dy,
\end{align*}
\\
with
\begin{align*}
f(x,y)=\frac{1}{2 \pi} \vert \Sigma^{(N)} \vert^{-\frac{1}{2}} \exp \left\{- \frac{(x,y)\Sigma^{(N)-1} (x,y)^t}{2} \right\}.
\end{align*}
\\
 where the matrix $\Sigma^{(N)}$ is the covariance matrix of $(\l_i,\l_j)$.
\\ 
 We will use the result of Equation \eqref{eq:Sigma_N}, which will be demonstrated in Appendix \ref{app2}, that is 
 \\
 \begin{align}
\Sigma^{(N)}=  \begin{pmatrix}
   1 +  \eta^\star\frac{A_N(i)}{\sqrt{N}}& \eta^\star\frac{B_N(i,j)}{\sqrt{N}} \\
  \eta^\star \frac{B_N(i,j)}{\sqrt{N}} & 1 + \eta^\star \frac{A_N(j)}{\sqrt{N}} 
\end{pmatrix},
\end{align} 
\\
where $A_N(i)=O_p(1)$, $A_N(j)=O_p(1)$ and $B_N(i,j)=O_p(1)$.
\\
We have
\\
\begin{align*}
f(x,y) &= \frac{1}{2\pi \vert \Sigma^{(N)} \vert^{-\frac{1}{2}}} \exp \left\{ -\frac{1}{2 \vert \Sigma^{(N)} \vert}\left[x^2(1+\frac{\eta^\star}{\sqrt{N}}A_N(j)) +
y^2(1+\frac{\eta^\star}{\sqrt{N}}A_N(i)) -2xy \frac{\eta^\star}{\sqrt{N}} B_N(i,j) \right] \right\} \\
&= \frac{1}{2\pi \vert \Sigma^{(N)} \vert^{-\frac{1}{2}}} \exp(-\frac{x^2}{2})\exp(-\frac{y^2}{2}) \exp \left\{ -\frac{x^2}{2} \left(\frac{1}{\vert \Sigma^{(N)}\vert} \left[1 + \frac{\eta^\star}{\sqrt{N}} A_N(j) \right]-1\right) \right.  \\
& \hspace{50mm}\left.-\frac{y^2}{2} \left(\frac{1}{\vert \Sigma^{(N)}\vert} \left[1 + \frac{\eta^\star}{\sqrt{N}} A_N(i) \right]-1\right) + \frac{1}{\vert \Sigma^{(N)}\vert} xy \frac{\eta^\star}{\sqrt{N}} B_N(i,j) \right\}.
\end{align*}
$$ $$ 
Using a first order Taylor development,

  $$\vert \Sigma^{(N)} \vert^{-1}= 1 - (A_N(i)+A_N(j))\frac{\eta^\star}{\sqrt{N}} + \alpha_N$$ and 

$$\vert \Sigma^{(N)} \vert^{-\frac{1}{2}}= 1 - \frac{1}{2}(A_N(i)+A_N(j))\frac{\eta^\star}{\sqrt{N}} + \beta_N,$$
where $\alpha_N=O_p(\frac{1}{N})$ and $\beta_N=O_p(\frac{1}{N})$.
\\ 
 More precisely,
\begin{align*}
 \alpha_N & =-(A_N(i) A_N(j) - B_N(i,j)^2)\frac{\eta^{\star2}}{N} \\
 &+ \frac{1}{2}\left(-(A_N(i)+A_N(j))\frac{\eta^\star}{\sqrt{N}}-(A_N(i) A_N(j) - B_N(i,j)^2) \frac{\eta^{\star2}}{N}\right)^2 \frac{1}{(1+\tilde{\alpha})^3},
 \end{align*} 
\\ 
 with $\vert\tilde{\alpha} \vert \leq \vert (A_N(i)+A_N(j))\frac{\eta^\star}{\sqrt{N}}+(A_N(i) A_N(j) - B_N(i,j)^2) \frac{\eta^{\star2}}{N} \vert .$
 \\
Similarly,
\\
\begin{align*}
 \beta_N & =-\frac{1}{2}(A_N(i) A_N(j) - B_N(i,j)^2)\frac{\eta^{\star2}}{N} \\
 & + \frac{1}{2}\left(-\frac{1}{2}(A_N(i)+A_N(j))\frac{\eta^\star}{\sqrt{N}}-\frac{1}{2}(A_N(i) A_N(j) - B_N(i,j)^2) \frac{\eta^{\star2}}{N}\right)^2 \frac{3}{4}\frac{1}{(1+\tilde{\beta})^\frac{5}{2}},
  \end{align*} 
 \\ 
   with $\vert\tilde{\beta} \vert \leq \vert \frac{1}{2} (A_N(i)+A_N(j))\frac{\eta^\star}{\sqrt{N}}+\frac{1}{2}(A_N(i) A_N(j) - B_N(i,j)^2) \frac{\eta^{\star2}}{N} \vert.$
 \\ Then,
\begin{align*}
f(x,y)& = \left(1 -\frac{1}{2} (A_N(i)+A_N(j))\frac{\eta^\star}{\sqrt{N}} + \beta_N \right) \phi(x)\phi(y) \\
& \times \exp \left\{ -\frac{x^2}{2}(-A_N(i)\frac{\eta^\star}{\sqrt{N}} +\gamma_N) -\frac{y^2}{2}(-A_N(j)\frac{\eta^\star}{\sqrt{N}} +\tilde{\gamma}_N) +  xy \left(\frac{\eta^\star}{\sqrt{N}} B_N(i,j) + \tilde{\tilde{\gamma}}_N \right) \right\}
\end{align*}
\\
where
$\gamma_N=-A_N(j)(A_N(i)+A_N(j))\frac{\eta^{\star2}}{N}+\alpha_N (1+ A_N(j)\frac{\eta^\star}{\sqrt{N}})=O_p\left(\frac{1}{N}\right)$,
\\
$\tilde{\gamma_N}=-A_N(i)(A_N(i)+A_N(j))\frac{\eta^{\star2}}{N}+\alpha_N (1+ A_N(i)\frac{\eta^\star}{\sqrt{N}})=O_p\left(\frac{1}{N}\right)$ and
\\
$\tilde{\tilde{\gamma_N}}=\frac{\eta^\star}{\sqrt{N}}B_N(i,j) \left( - (A_N(i)+A_N(j))\frac{\eta^\star}{\sqrt{N}} + \alpha_N \right)=O_p\left(\frac{1}{N}\right)$
\\
A Taylor development of the exponential function leads to
\\
\begin{align*}
f(x,y)& =\left(1 - \frac{1}{2}(A_N(i)+A_N(j))\frac{\eta^\star}{\sqrt{N}} + \beta_N \right) \phi(x)\phi(y) \\
 &  \times \left[ 1+\frac{x^2}{2}\frac{\eta^\star}{\sqrt{N}} A_N(i) +\frac{y^2}{2}\frac{\eta^\star}{\sqrt{N}} A_N(j) + xy \frac{\eta^\star}{\sqrt{N}} B_N(i,j) +\nu_N(x)  \right]
\end{align*}
\\
with
\\
\begin{align*}
 \nu_N(x)& =-\frac{x^2}{2} \gamma_N -\frac{y^2}{2} \tilde{\gamma}_N +xy \tilde{\tilde{\gamma}}_N \\
 &+ \frac{1}{2} \left(\frac{x^2}{2}\frac{\eta^\star}{\sqrt{N}} A_N(i) +\frac{y^2}{2}\frac{\eta^\star}{\sqrt{N}} A_N(j) + xy \frac{\eta^\star}{\sqrt{N}} B_N(i,j) - \frac{x^2}{2} \gamma_N -\frac{y^2}{2} \tilde{\gamma}_N +xy \tilde{\tilde{\gamma}}_N \right)^2 \exp{\tilde{u}}
 \end{align*}
\\ 
 where $\vert \tilde{u} \vert \leq \vert \frac{x^2}{2}\frac{\eta^\star}{\sqrt{N}} A_N(i) +\frac{y^2}{2}\frac{\eta^\star}{\sqrt{N}} A_N(j) + xy \frac{\eta^\star}{\sqrt{N}} B_N(i,j) - \frac{x^2}{2} \gamma_N -\frac{y^2}{2} \tilde{\gamma_N} +xy \tilde{\tilde{\gamma_N}} \vert.$
\\
Then, 
\\
\begin{align*}
\int_{t}^{\infty} \int_{t}^{\infty} f(x,y) dx dy & =  \left(1 -\frac{1}{2} (A_N(i)+A_N(j))\frac{\eta^\star}{\sqrt{N}} + \beta_N \right) \\ & \left[K^2+\frac{1}{2}\frac{\eta^\star}{\sqrt{N}} (A_N(j) + A_N(i))K(K +t\phi(t)) + B_N(i,j) \frac{\eta^\star}{\sqrt{N}} \phi(t)^2  \right] + \mu_N \\ 
& = K^2 + \frac{1}{2}(A_N(i)+A_N(j))\frac{\eta^\star}{\sqrt{N}} Kt\phi(t) + B_N(i,j) \frac{\eta^\star}{\sqrt{N}} \phi(t)^2 +  \mu_N'
\end{align*}
\\
where $\mu_N=\left(1 -\frac{1}{2} (A_N(i)+A_N(j) +\beta_N) \frac{\eta^\star}{\sqrt{N}} \right)\int_{t}^{\infty} \int_{t}^{\infty} \phi(x)\phi(y) \nu_N(x) dx dy$ 
\begin{align*}
\textrm{and }\mu_N'&=\mu_N + \beta_N \left(K^2+\frac{1}{2}\frac{\eta^\star}{\sqrt{N}} (A_N(j) + A_N(i))K(K +t\phi(t)) + B_N(i,j) \frac{\eta^\star}{\sqrt{N}} \phi(t)^2\right) \\
& - \frac{1}{2} (A_N(i)+A_N(j))\frac{\eta^{\star2}}{N} B_N(i,j) \phi(t)^2.
\end{align*}
\\
This remainder and its order will be carefully studied in Section \ref{sec:proof_restes}.
\\
 Similarly, we can compute  $\mathbb{P}(\Y_i=\Y_j=0 \vert \Z)$ and  $\mathbb{P}(\Y_i \neq \Y_j \vert \Z)$:
\\
\begin{align*}
\int_{-\infty}^{t} \int_{-\infty}^{t} f(x,y) dx dy 
 = (1-K)^2 - \frac{1}{2}(A_N(i)+A_N(j))\frac{\eta^\star}{\sqrt{N}} (1-K)t\phi(t) + B_N(i,j) \frac{\eta^\star}{\sqrt{N}} \phi(t)^2 +  \tilde{\mu}_N
\end{align*}
 
 \begin{align*}
\int_{-\infty}^{t} \int_{t}^{\infty} f(x,y) dx dy + \int_{t}^{\infty} \int_{-\infty}^{t} f(x,y) dx dy =2 K(1-K) & + (A_N(i)+A_N(j))\frac{\eta^\star}{\sqrt{N}}(1-2K)t\phi(t) \\
& -2 B_N(i,j) \frac{\eta^\star}{\sqrt{N}} \phi(t)^2 +  \tilde{\tilde{\mu}}_N.
\end{align*}

$$ $$
Replacing these terms in the expression of the numerator of $\mathbb{E}[\W_i \W_j \vert \Z, \epsilon_i=\epsilon_j=1] $ given in equation (\ref{E_Wi_Wj}) leads to:

\begin{align}
\frac{\eta^\star}{\sqrt{N}} B_N(i,j) \phi(t)^2 \frac{(1-P)}{P(1-K)^2} + r_N,
\end{align}
\\
 where $r_N$ is a linear combination of $\mu_N'$, $\tilde{\mu}_N$ and $\tilde{\tilde{\mu}}_N$.
  \\
Since there is no constant term in this numerator, we only need the development of order 0 of the denominator of $\mathbb{E}[\W_i \W_j \vert \Z, \epsilon_i=\epsilon_j=1] $ to obtain the first order approximation of $\mathbb{E}[\W_i \W_j \vert \Z, \epsilon_i=\epsilon_j=1] $.
\\
We obtain that the denominator can be written as 
\\
$$\frac{K^2}{P^2} +\tilde{r}_N,$$
\\
where $\tilde{r}_N$ is the sum of a term of order $\frac{1}{\sqrt{N}}$ and a linear combination of $\mu_N'$, $\tilde{\mu}_N$ and $\tilde{\tilde{\mu}}_N$.
Thus, we obtain that 
\\
\begin{align}
\mathbb{E}[\W_i \W_j \vert \Z, \epsilon_i=\epsilon_j=1]& =\frac{\frac{\eta^\star}{\sqrt{N}} B_N(i,j) \phi(t)^2 \frac{(1-P)}{P(1-K)^2} + r_N }{\frac{K^2}{P^2} + \tilde{r}_N} \\
&= \eta^\star \G_N(i,j) \phi(t)^2 \frac{P(1-P)}{K^2(1-K)^2} + R_N(i,j)
\end{align}
\\
where 
\\
\begin{equation}
R_N(i,j)=\left(\frac{\eta^\star}{\sqrt{N}} B_N(i,j) \phi(t)^2 \frac{(1-P)}{P(1-K)^2} + r_N \right)\tilde{r}_N + \frac{K^2}{P^2} r_N.
\label{eq:R_N}
\end{equation}
\subsection{Proof of Theorem \ref{th:consistency}}
\label{sec:proof_consistency}

\subsubsection{Properties of $\Z$}

In the following proofs, we will use several properties of the matrix $\Z$, which are stated in Proposition \ref{prop:prop_Z_EJS}.

\begin{prop}
Uniformly in k,
\begin{enumerate}[label=(\arabic*)]
\item \label{prop1:1} $\mathbb{E}[\Z_{1,k}\Z_{2,k}]=-\frac{1}{n-1}$. 
\item  \label{prop1:2} $\mathbb{E}[\Z_{1,k}^p]=O(1)$, for all $p$.
\item \label{prop1:3} $\mathbb{E}[\Z_{1,k}^2\Z_{2,k}^2]= 1 + o(1)$.
\item \label{prop1:4} $\mathbb{E}[\Z_{1,k}^3\Z_{2,k}]=  O\left(\frac{1}{n}\right)$.
\item \label{prop1:5}$ \mathbb{E}[\Z_{1,k}^2 \Z_{2,k} \Z_{3,k}]  = O\left(\frac{1}{n}\right)$.
\item \label{prop1:6}$ \mathbb{E}[\Z_{1,k} \Z_{2,k} \Z_{3,k} \Z_{4,k}] = O\left(\frac{1}{n^2}\right)$.
\item \label{prop1:7}$\mathbb{E}[\Z_{1,k}^5\Z_{2,k}]=O\left(\frac{1}{n}\right)$.
\item \label{prop1:8}$\mathbb{E}[\Z_{1,k}^3\Z_{2,k}^3]=O(1)$.
\item \label{prop1:9}$\mathbb{E}[\Z_{1,k}^4\Z_{2,k}^2]=O(1)$.
\item \label{prop1:10}$\mathbb{E}[\Z_{1,k}^4\Z_{2,k}\Z_{3,k}]=O(\frac{1}{ n})$.
\item \label{prop1:11}$\mathbb{E}[\Z_{1,k}^3\Z_{2,k}^2\Z_{3,k}]=O(\frac{1}{n})$.
\item \label{prop1:12}$\mathbb{E}[\Z_{1,k}^3\Z_{2,k}\Z_{3,k}\Z_{4,k}]=O(\frac{1}{n^2})$.
\end{enumerate}
\label{prop:prop_Z_EJS}
\end{prop}
The proof of Proposition \ref{prop:prop_Z_EJS} is given in Appendix \ref{app4}.
%
%
%
%
%
%

\subsubsection{Proof of Lemma \ref{lemma:denominator}}
\label{sec:denominator}

Let us prove that, when $n$ and $N$ go to infinity and $n/N$ goes to $a$,  $$\frac{1}{n}\underset{i \neq j}{\sum}  \G_N(i,j)^2 \overset{P}{\rightarrow} a, $$ 
\\
where $\overset{P}{\rightarrow} $ denotes the convergence in probability.
\\
\begin{align*}
\G_N(i,j)^2= \frac{1}{N^2} \underset{k=1}{\overset{N}{\sum}}\Z_{i,k}^2\Z_{j,k}^2 + \frac{1}{N^2} \underset{k \neq l}{\sum} \Z_{i,k}\Z_{j,k} \Z_{i,l} \Z_{j,l}
\end{align*}
\\
Since $\Z_{i,k}$ and $\Z_{j,l}$ are independent for any $i$ and $j$ when $k \neq l$, we will always consider separately the cases where $k=l$ from the cases where $k\neq l$.
\\
Indeed, let us show that 

\begin{equation}
\frac{1}{n}\frac{1}{N^2} \underset{i \neq j}{\sum} \underset{k=1}{\overset{N}{\sum}}\Z_{i,k}^2\Z_{j,k}^2 \overset{P}{\rightarrow} a
\label{eq:conv_a}
\end{equation}
\\
and
\begin{equation}
 \frac{1}{n}\frac{1}{N^2} \underset{i \neq j}{\sum} \underset{k \neq l}{\sum} \Z_{i,k}\Z_{j,k} \Z_{i,l} \Z_{j,l}\overset{P}{\rightarrow} 0.
\label{eq:conv_0}
\end{equation} 
 
%
Note that

\begin{align*}
 \mathbb{E}[\frac{1}{n}\frac{1}{N^2} \underset{i \neq j}{\sum} \underset{k=1}{\overset{N}{\sum}}\Z_{i,k}^2\Z_{j,k}^2] & = \frac{1}{n}\frac{1}{N^2} \underset{i \neq j}{\sum} \underset{k=1}{\overset{N}{\sum}}\mathbb{E}[\Z_{i,k}^2\Z_{j,k}^2] \\
& = \frac{n-1}{N}(1+o(1)) \textrm{ by  \ref{prop1:3} of Proposition \ref{prop:prop_Z_EJS}} \\
& = a +o(1)
\end{align*}
\\
Moreover, 
\begin{align}
 \textrm{Var} (\frac{1}{n}\frac{1}{N^2} \underset{i \neq j}{\sum} \underset{k=1}{\overset{N}{\sum}}\Z_{i,k}^2\Z_{j,k}^2)& = \frac{1}{n^2}\frac{1}{N^4}\underset{k=1}{\overset{N}{\sum}} \underset{i_1 \neq j_1}{\sum}\underset{i_2 \neq j_2}{\sum} \mathbb{E}[\Z_{i_1,k}^2\Z_{j_1,k}^2 \Z_{i_2,k}^2\Z_{j_2,k}^2] - \frac{1}{n^2}\frac{1}{N^4}\underset{k=1}{\overset{N}{\sum}} \left(\underset{i \neq j}{\sum} \mathbb{E}[\Z_{i,k}^2\Z_{j,k}^2]\right)^2 
 \label{eq:var_lemma1}
\end{align}
\\
The second term of \eqref{eq:var_lemma1} can be rewritten as:
\\
\begin{align*}
\frac{1}{n^2}\frac{1}{N^4}\underset{k=1}{\overset{N}{\sum}} \left(\underset{i \neq j}{\sum} \mathbb{E}[\Z_{i,k}^2\Z_{j,k}^2]\right)^2 & = \frac{N n^2 (n-1)^2}{n^2 N^4} (1+o(1)) \textrm{ by \ref{prop1:3} of Proposition \ref{prop:prop_Z_EJS}}\\
& = O \left(\frac{1}{n}\right)
\end{align*}
\\
\begin{align*}
\underset{i_1 \neq j_1}{\sum}\underset{i_2 \neq j_2}{\sum} \mathbb{E}[\Z_{i_1,k}^2\Z_{j_1,k}^2 \Z_{i_2,k}^2\Z_{j_2,k}^2] \leq  \mathbb{E}[\underset{i_1, j_1,i_2,j_2}{\sum}\Z_{i_1,k}^2\Z_{j_1,k}^2 \Z_{i_2,k}^2\Z_{j_2,k}^2]= \mathbb{E}\left[\underset{i=1}{\overset{n}{\sum}} \Z_{i,k}^2\right]^4=n^4
\end{align*}
\\
This last equality comes from the definition of $\Z$ as a centered and normalized variable given in Equation \eqref{eq:normalization_1}, which implies that for all $k$,

$$ \underset{i=1}{\overset{n}{\sum}} \Z_{i,k}^2 =n.$$
\\
Then,

$$\frac{1}{n^2}\frac{1}{N^4}\underset{k=1}{\overset{N}{\sum}} \underset{i_1 \neq j_1}{\sum}\underset{i_2 \neq j_2}{\sum} \mathbb{E}[\Z_{i_1,k}^2\Z_{j_1,k}^2 \Z_{i_2,k}^2\Z_{j_2,k}^2] \leq \frac{n^4 N}{n^2 N^4} = O\left(\frac{1}{n}\right).$$
\\
This proves \eqref{eq:conv_a}.


\begin{align*}
\mathbb{E}[\frac{1}{n}\frac{1}{N^2} \underset{i \neq j}{\sum} \underset{k \neq l}{\sum} \Z_{i,k}\Z_{j,k} \Z_{i,l} \Z_{j,l}] & =\frac{1}{n}\frac{1}{N^2} \underset{i \neq j}{\sum} \underset{k \neq l}{\sum} \mathbb{E}[ \Z_{i,k}\Z_{j,k}]\mathbb{E}[\Z_{i,l} \Z_{j,l}] \\
& = \frac{n (n-1) N(N-1)}{n N^2 (n-1)^2} \textrm{ by  \ref{prop1:1} of Proposition \ref{prop:prop_Z_EJS}} \\
& = O\left(\frac{1}{n}\right)
\end{align*}
%
 
 \begin{align*}
\textrm{Var}(\frac{1}{n}\frac{1}{N^2} \underset{i \neq j}{\sum} \underset{k \neq l}{\sum} \Z_{i,k}\Z_{j,k} \Z_{i,l} \Z_{j,l}) 
& = \frac{1}{n^2}\frac{1}{N^4}  \underset{k \neq l}{\sum} \underset{i_1 \neq j_1}{\sum}\underset{i_2 \neq j_2}{\sum} \mathbb{E}[\Z_{i_1,k}\Z_{i_2,k}\Z_{j_1,k}\Z_{j_2,k}]\mathbb{E}[\Z_{i_1,l}\Z_{i_2,l}\Z_{j_1,l}\Z_{j_2,l}] \\
& - \frac{1}{n^2}\frac{1}{N^4}  \underset{k \neq l}{\sum} \left( \underset{i \neq j}{\sum} \mathbb{E}[\Z_{i,k}\Z_{j,k}]\mathbb{E}[ \Z_{i,l} \Z_{j,l}] \right)^2
\end{align*}
 
 \begin{align*}
 \frac{1}{n^2}\frac{1}{N^4}  \underset{k \neq l}{\sum} \left( \underset{i \neq j}{\sum} \mathbb{E}[\Z_{i,k}\Z_{j,k}]\mathbb{E}[ \Z_{i,l} \Z_{j,l}] \right)^2 & = \frac{N(N-1)n^2 (n-1)^2}{n^2 N^4 (n-1)^4}  \textrm{ by \ref{prop1:1} of Proposition \ref{prop:prop_Z_EJS}} \\
 & = O\left(\frac{1}{n^4}\right)
  \end{align*}
 \\
In the first term, $\{i_1,i_2,j_1,j_2 \}$ can be of cardinal 2, 3 or 4 and counting the number of combinations gives the expression: 
 \begin{align*}
 \underset{i_1 \neq j_1}{\sum}\underset{i_2 \neq j_2}{\sum} & \mathbb{E}[\Z_{i_1,k}\Z_{i_2,k}\Z_{j_1,k}\Z_{j_2,k}]\mathbb{E}[\Z_{i_1,l}\Z_{i_2,l}\Z_{j_1,l}\Z_{j_2,l}] = 2 \underset{i\neq j}{\sum} \mathbb{E}[\Z_{i,k}^2 \Z_{j,k}^2]\mathbb{E}[\Z_{i,l}^2 \Z_{j,l}^2] \\
 & + 4  \underset{i\neq j_1 \neq j_2} {\sum} \mathbb{E}[\Z_{i,k}^2 \Z_{j_1,k} \Z_{j_2,k}] \mathbb{E}[\Z_{i,l}^2 \Z_{j_1,l} \Z_{j_2,l}]  \\
 &+ \underset{i_1 \neq i_2 \neq j_1 \neq j_2} {\sum}  \mathbb{E}[\Z_{i_1,k}\Z_{i_2,k}\Z_{j_1,k}\Z_{j_2,k}]\mathbb{E}[\Z_{i_1,l}\Z_{i_2,l}\Z_{j_1,l}\Z_{j_2,l}] \\
 & = 2n(n-1) (1+ o(1)) + 4 \frac{n(n-1)(n-2)}{n} o(1) + \frac{n(n-1)(n-2)(n-3)}{n^2}o(1)  =O(n^2)
 \end{align*}
\\ 
 This was obtained by using \ref{prop1:3},\ref{prop1:5} and \ref{prop1:6} of Proposition  \ref{prop:prop_Z_EJS}.
\\ 
 Finally, 
 
 $$ Var\left(\frac{1}{n}\frac{1}{N^2} \underset{i \neq j}{\sum} \underset{k \neq l}{\sum} \Z_{i,k}\Z_{j,k} \Z_{i,l} \Z_{j,l}\right) = O\left(\frac{1}{n^2}\right).$$
\\ 
 This completes the proof of \eqref{eq:conv_0}.
 
%
%
%
%
%
%
 
  \subsubsection{Proof of Lemma \ref{lemma5}}

Note that
\begin{align*}
\mathbb{P}(E_N^c) & \leq n \underset{i}{\textrm{ sup }} \mathbb{P}\left(\vert \overset{N}{\underset{k=1}{\sum}} (\Z_{i,k}^2-1) \vert \geq N \epsilon_N\right) + n (n-1) \underset{i \neq j}{\textrm{ sup }} \mathbb{P}\left(\vert \overset{N}{\underset{k=1}{\sum}} \Z_{i,k} \Z_{j,k} \vert \geq N \epsilon_N\right) \\
& = n \mathbb{P}\left(\vert \overset{N}{\underset{k=1}{\sum}} (\Z_{1,k}^2-1) \vert \geq N \epsilon_N\right) + n(n-1) \mathbb{P}\left(\vert \overset{N}{\underset{k=1}{\sum}} \Z_{1,k} \Z_{2,k} \vert \geq N \epsilon_N\right).
\end{align*}
\\
Let $\delta$ be a positive real number such that $\sqrt{\delta}/2c \in V_0$ and $\delta < \frac{\delta_{min}}{4}$, where $V_0$ and $\delta_{min}$ are defined in Assumptions \ref{HYP:HYP1} and \ref{hyp3} respectively.

\begin{align*}
\mathbb{P}\left(\vert \overset{N}{\underset{k=1}{\sum}} (\Z_{i,k}^2-1) \vert \geq N \epsilon_N\right) \leq \mathbb{P}\left( \exists k, s_k^2 \leq \delta \right) + \mathbb{P}\left(\vert \overset{N}{\underset{k=1}{\sum}} (A_{i,k}-\bar{A}_k)^2 -s_k^2) \vert \geq N \delta \epsilon_N \right) 
\end{align*}
\\
Note also that
\begin{align*}
\left\{\exists k, s_k^2 \leq \delta \right\} = \overset{N}{\underset{k=1}{\bigcup}} \left\{ \overset{n}{\underset{i=1}{\sum}} (A_{i,k}-\bar{A}_k)^2 \leq n\delta \right\} = \overset{N}{\underset{k=1}{\bigcup}}\left\{ \overset{N}{\underset{i=1}{\sum}} (A_{i,k}-m_k +m_k - \bar{A}_k)^2 \leq n\delta \right\}
\end{align*}

where $m_k=\mathbb{E}[A_{i,k}]$. \\

Observe that
\begin{align}
\left\{ \overset{n}{\underset{i=1}{\sum}} (A_{i,k}-m_k +m_k - \bar{A}_k)^2 \leq n\delta \right\} \subset \left\{ \vert \bar{A}_k -m_k \vert \geq \sqrt{\delta} \right\} \cup \left\{ \overset{n}{\underset{i=1}{\sum}} (A_{i,k}-m_k)^2 \leq 4n\delta \right\}.
\label{eq:inclu_union}
\end{align}

%
%
%
%
%

Let us show that 

\begin{equation}
\mathbb{P}(\vert \bar{A}_k -m_k \vert \geq \sqrt{\delta}) \leq 2C \exp \left\{ -\frac{n\delta}{4d} \right\}.
\label{eq:control1}
\end{equation}

%
%
%
%
%

$$\mathbb{P}(\vert \bar{A}_k -m_k \vert \geq \sqrt{\delta}) = \mathbb{P}( \bar{A}_k -m_k  \geq \sqrt{\delta}) + \mathbb{P}( \bar{A}_k -m_k  \leq -\sqrt{\delta})$$

By Chernoff inequality, for all $\lambda \geq 0$,

\begin{align*}
 \mathbb{P}(n( \bar{A}_k -m_k)  \geq n \sqrt{\delta}) &\leq \exp\left\{ -n \sqrt{\delta}\lambda + \log \left( \mathbb{E}[ \exp(n (\bar{A}_k -m_k))] \right)\right\} \\
 &= \exp\left\{ -n \sqrt{\delta}\lambda +  n \log \left( \mathbb{E}[ \exp(A_{i,k} -m_k)] \right)\right\}
\end{align*}

Then, by Assumption \ref{hyp2}, for all positive values of $ \lambda$ in $V_0$,

\begin{align}
 \mathbb{P}(n (\bar{A}_k -m_k)  \geq n \sqrt{\delta}) \leq C \exp\left\{ -n \sqrt{\delta}\lambda +  n d \lambda^2 \right\}.
 \label{eq:chernoff1}
\end{align}

The right term of \eqref{eq:chernoff1} is maximum when 

$$\lambda= \frac{\sqrt \delta}{2d},$$

which implies that 

$$\mathbb{P}( \bar{A}_k -m_k \geq \sqrt{\delta}) \leq C \exp \left\{ -\frac{n\delta}{4d} \right\}.$$ 

Similarly, for all negative values of $ \lambda$ in $V_0$,
\begin{align}
 \mathbb{P}(n (\bar{A}_k -m_k)  \leq - n \sqrt{\delta}) \leq C \exp\left\{ n \sqrt{\delta}\lambda +  n d \lambda^2 \right\}.
  \label{eq:chernoff2}
\end{align}

The right term of \eqref{eq:chernoff2} is maximum when 

$$\lambda= -\frac{\sqrt \delta}{2d},$$

which implies that 

$$\mathbb{P}( \bar{A}_k -m_k \geq \sqrt{\delta}) \leq C \exp \left\{ -\frac{n\delta}{4d} \right\},$$ 

which proves \eqref{eq:control1}.


$$\mathbb{P}(\overset{n}{\underset{i=1}{\sum}} (A_{i,k}-m_k)^2 \leq 4n\delta ) \leq \mathbb{P}(\overset{n}{\underset{i=1}{\sum}}[ (A_{i,k}-m_k)^2 -\sigma_k^2] \leq n(4\delta - \delta_{min}) )$$

Since $4\delta - \delta_{min}<0$ by assumption on $\delta$, we apply again Chernoff inequality, which gives us that:

\begin{align*}
\mathbb{P}(\overset{n}{\underset{i=1}{\sum}}[ (A_{i,k}-m_k)^2 -\sigma_k^2] & \leq n(4\delta - \delta_{min}) ) \leq C \exp\left\{ -n \frac{(4\delta-\delta_{min})^2}{2d} \right\} 
\end{align*}


This result, combined with \eqref{eq:control1}, proves that

\begin{align}
\mathbb{P}\left( \exists k, s_k^2 \leq \delta \right)& \leq 2NC \exp \left\{ -\frac{n\delta}{4d} \right\} + NC \exp\left\{ -n \frac{(4\delta-\delta_{min})^2}{2d} \right\} 
\label{eq:control_sk}
\end{align}
%
 
Notice that 
 \begin{align*} 
& \left\{ \left\vert\overset{N}{\underset{k=1}{\sum}} (A_{i,k}-\bar{A}_k)^2 -s_k^2 \right\vert\geq N\delta \epsilon_N \right\}= \left\{ \frac{1}{n}\left\vert \overset{N}{\underset{k=1}{\sum}} \overset{n}{\underset{l=1}{\sum}} (A_{i,k}-\bar{A}_k)^2 -(A_{l,k}-\bar{A}_k)^2 \right\vert \geq N\delta \epsilon_N \right\}  \\
& \subset \left\{ \left\vert\overset{N}{\underset{k=1}{\sum}} (A_{i,k}-m_k)^2 -\sigma_k^2 \right\vert\geq \frac{N\delta \epsilon_N}{4} \right\} \cup \left\{ \left\vert \overset{N}{\underset{k=1}{\sum}} \overset{n}{\underset{l=1}{\sum}} (A_{l,k}-m_k)^2 -\sigma_k^2 \right\vert \geq \frac{n N\delta \epsilon_N}{4} \right\}  \\
& \cup \left\{ \left\vert\overset{N}{\underset{k=1}{\sum}} (A_{i,k}-m_k)(m_k-\bar{A}_k) \right\vert\geq \frac{N\delta \epsilon_N}{8} \right\} \cup \left\{ \left\vert \overset{N}{\underset{k=1}{\sum}} \overset{n}{\underset{l=1}{\sum}} (A_{l,k}-m_k)(m_k-\bar{A}_k) \right\vert \geq \frac{n N\delta \epsilon_N}{8} \right\}
\end{align*} 
 
 Using Chernoff inequality and Assumption \ref{hyp1}, we can prove that

 \begin{align*}
 \mathbb{P}\left(\left\vert\overset{N}{\underset{k=1}{\sum}} (A_{i,k}-m_k)^2 -\sigma_k^2 \right\vert\geq \frac{N\delta \epsilon_N}{4} \right) \leq 2C \exp \left\{-\frac{N\delta^2\epsilon_N^2}{64d} \right\}
 \end{align*}

 and

 \begin{align*}
 \mathbb{P}\left(\left\vert \overset{N}{\underset{k=1}{\sum}} \overset{n}{\underset{l=1}{\sum}} (A_{l,k}-m_k)^2 -\sigma_k^2 \right\vert\geq \frac{n N\delta \epsilon_N}{4} \right) \leq 2C \exp \left\{-\frac{N n \delta^2\epsilon_N^2}{64d} \right\}
 \end{align*}

 Moreover,
 
 \begin{align*}
\mathbb{P}\left( \left\vert\overset{n}{\underset{k=1}{\sum}} (A_{i,k} - m_k)( m_k -\bar{A}_k) \right\vert \geq \frac{N \delta \epsilon_N}{4} \right) &
\leq \mathbb{P}\left( \overset{n}{\underset{k=1}{\sum}} (A_{i,k} - m_k)^2  \geq N n \frac{\delta \epsilon_N}{8}\right)
\\
& + \mathbb{P}\left(\left\vert \overset{n}{\underset{k=1}{\sum}} \underset{l \neq i}{\sum}(A_{i,k} - m_k)( m_k -A_{l,k}) \right\vert \geq n N \frac{\delta \epsilon_N}{8}\right)
\end{align*}

Using Chernoff inequality and Assumption \ref{hyp4}, we obtain that
 $$ \mathbb{P}\left(\left\vert \overset{n}{\underset{k=1}{\sum}} \underset{l \neq i}{\sum}(A_{i,k} - m_k)( m_k -A_{l,k}) \right\vert \geq n N \frac{\delta \epsilon_N}{8}\right) \leq 2 C \exp \left\{ - \frac{n N\delta^2 \epsilon_N^2}{256d} \right\}$$

and with Assumption \ref{hyp1} we have
\\
\begin{align}
 \mathbb{P}\left(\overset{n}{\underset{k=1}{\sum}} (A_{i,k} - m_k)^2 \geq N n \frac{\delta \epsilon_N}{8}\right)
 & \leq  C \exp \left\{ - \frac{n^2 N\delta^2 \epsilon_N^2 }{256 d }+ \frac{n N \delta \delta_{max} \epsilon_N}{16d} - \frac{N \delta_{max}}{4d} \right\}, 
 \label{eq:last_control}
 \end{align}
 
with  $n^2 N \epsilon_N^2 = a^2 	N^{2 +2 \gamma}$ and $ n N \epsilon_N = a N^{\frac{3}{2}+\gamma}$ where $\gamma >0$, which implies that the main term in the exponential is 
$- \frac{n^2 N\delta^2 \epsilon_N^2 }{256 d }$.
 
 Similarly, we can show that 
 
 \begin{align*}
 \mathbb{P}\left(\left\vert \overset{N}{\underset{k=1}{\sum}} \overset{n}{\underset{l=1}{\sum}} (A_{l,k}-m_k)(m_k-\bar{A}_k) \right\vert \geq \frac{n N\delta \epsilon_N}{8} \right)
 & \leq 2C \exp \left\{ - \frac{n^2 N \delta^2 \epsilon_N^2}{256d} \right\} \\
 &+ C \exp \left\{ - \frac{n^3 N\delta^2 \epsilon_N^2}{256 d } 
 + \frac{n^2 N \delta \delta_{max} \epsilon_N}{16d} - \frac{Nn \delta_{max}}{4d} \right\}.
 \end{align*}

 This concludes the proof that for all values of $q$,

 $$  \mathbb{P}\left(\vert \overset{n}{\underset{k=1}{\sum}} (\Z_{i,k}^2-1) \vert \geq N \epsilon_N\right) = O \left(\frac{1}{N^q}\right).$$

 We use similar techniques to otain an upper bound for $\mathbb{P}\left(\left\vert \overset{n}{\underset{k=1}{\sum}} \Z_{i,k} \Z_{j,k} \right\vert \geq N \epsilon_N\right)$.
\begin{align*}
\mathbb{P}\left(\left\vert \overset{n}{\underset{k=1}{\sum}} \Z_{i,k} \Z_{j,k} \right\vert \geq N \epsilon_N\right)&  = \mathbb{P}\left(\left\vert \overset{n}{\underset{k=1}{\sum}} \frac{(A_{i,k} - \bar{A}_k)(A_{j,k} - \bar{A}_k)}{s_k^2} \right\vert \geq N \epsilon_N\right) \\
& \leq \mathbb{P}(\exists k, s_k^2 \leq \delta) + \mathbb{P}\left(\left\vert \overset{n}{\underset{k=1}{\sum}} (A_{i,k} - \bar{A}_k)(A_{j,k} - \bar{A}_k) \right\vert \geq N \delta \epsilon_N\right)
\end{align*}

Since we have already proved \eqref{eq:control_sk} and \eqref{eq:last_control}, we will conclude the proof by showing that

\begin{align}
\mathbb{P}\left(\left\vert \overset{n}{\underset{k=1}{\sum}} (A_{i,k} - m_k)(A_{j,k} - m_k) \right\vert \geq N \frac{\delta \epsilon_N}{4}\right)\leq 2C \exp \left\{ - \frac{N\delta \epsilon_N}{64 d} \right \},
\label{eq:control_3}
\end{align}

and

\begin{align}
\mathbb{P}\left( \overset{n}{\underset{k=1}{\sum}} (\bar{A}_k - m_k)^2 \geq N \frac{\delta \epsilon_N}{4}\right) \leq  N^2 C \exp \left\{ - \frac{N\delta \epsilon_N}{16 d} \right \}.
\label{eq:control_5}
\end{align}

%
%
%
%


\eqref{eq:control_3} is obtained using Assumption \ref{hyp4} and Chernoff inequality.

\begin{align*}
\mathbb{P}\left( \overset{n}{\underset{k=1}{\sum}} (\bar{A}_k - m_k)^2  \geq N \frac{\delta \epsilon_N}{4}\right) & \leq \mathbb{P}\left(\underset{k}{\textrm{ sup }} (m_k-\bar{A}_k)^2 \geq \frac{\delta \epsilon_N}{4} \right) \\
& \leq N \underset{k}{\textrm{ sup }} \mathbb{P}\left( (m_k-\bar{A}_k)^2 \geq \frac{\delta \epsilon_N}{4} \right) \\
& \leq  N^2 C \exp \left\{ - \frac{N\delta \epsilon_N}{16 d} \right \},
\end{align*}

which proves \eqref{eq:control_5} and achieves the proof of Lemma \ref{lemma5}.

\subsubsection{Proof of Lemma \ref{lemma3}}
\label{sec:proof_restes}
According to the results of Section \ref{sec:denominator}, we have
\\
\begin{align*}
  \frac{1}{n}\underset{i \neq j}{\sum} \mathbb{E}[\W_i \W_j \vert \Z, \epsilon_i=\epsilon_j=1 ]\G_N(i,j)\mathbbm{1}_{E_N} & = \frac{1}{n}\underset{i \neq j}{\sum} (c\eta^\star \G_N(i,j) + R_N(i,j) ) \G_N(i,j) \mathbbm{1}_{E_N} \\
  & = ac \eta^\star + \frac{1}{n}\underset{i \neq j}{\sum}  R_N(i,j) \G_N(i,j)\mathbbm{1}_{E_N} + o_p(1)
\end{align*}
 \\ 
   Thus, we just need to prove that $\underset{i\neq j}{\sum} \G_N(i,j)R_N(i,j)\mathbbm{1}_{E_N}=o_p(1). $
 \\
%
%
We shall see that $R_N(i,j)\mathbbm{1}_{E_N}$ may be upper bounded by a finite sum of terms of the form

\begin{equation}
\vert \G_N(i,j)\vert^{k_1}\vert \G_N(i,i)-1 \vert^{k_2}\vert \G_N(j,j)-1 \vert^{k_3},
\label{eq:proof_RN}
\end{equation}
with k in $\llbracket 2, 22 \rrbracket$ and $k_1+k_2+k_3=k$.

Thus,  $\frac{1}{n}\underset{i \neq j}{\sum} R_N(i,j) \G_N(i,j)\mathbbm{1}_{E_N}$ is upper bounded by a finite sum of terms of the form 
$$ \frac{1}{n}\underset{i \neq j}{\sum} \vert \G_N(i,j)\vert^{k_1+1}\vert \G_N(i,i)-1 \vert^{k_2}\vert \G_N(j,j)-1 \vert^{k_3}.$$

But 
\begin{align*}
\frac{1}{n}\underset{i \neq j}{\sum} \vert \G_N(i,j)\vert^{k_1+1}\vert \G_N(i,i)-1 \vert^{k_2}\vert \G_N(j,j)-1 \vert^{k_3} \mathbbm{1}_{E_N} & \leq \epsilon_N^{k_1+k_2+k_3+1}  \frac{n(n-1)}{n} 
\\
&= O\left( \frac{1}{N^{\frac{1}{2}-3 \gamma}} \right) \\
&=o(1),
\end{align*}
since $k_1+k_2+k_3 +1 \geq 3$ and $ \gamma < 1/10$.
\\
This achieves the proof of Lemma \ref{lemma3}.
 \\
 Let us explain why Equation \eqref{eq:proof_RN} holds.
\\
%
%
%
%
%
%
  We need to evaluate $\vert R_N(i,j) \mathbbm{1}_{E_N} \vert$. Then, let us look at the previous remainders which compose $R_N(i,j)$, and we will provide upper bounds when $E_N$ holds.

  \begin{align*}
  \vert \alpha_N \vert = \vert A_N(i) A_N(j) - B_N(i,j)^2 \vert \frac{\eta^{\star2}}{N} + \frac{1}{2}\vert(A_N(i)+A_N(j))\frac{\eta^\star}{\sqrt{N}}+(A_N(i) A_N(j) - B_N(i,j)^2) \frac{\eta^{\star2}}{N}\vert^2 \frac{1}{\vert 1+\tilde{\alpha}\vert^3},
  \end{align*}
  
  with $\vert\tilde{\alpha} \vert \leq \vert (A_N(i)+A_N(j))\frac{\eta^\star}{\sqrt{N}}+(A_N(i) A_N(j) - B_N(i,j)^2) \frac{\eta^{\star2}}{N} \vert \leq 2\epsilon_N \eta^\star + 2 \epsilon_N^2 \eta^{\star 2}$.
  
  Similarly,
   \begin{align*}
  \vert \beta_N \vert & = \frac{1}{2}\vert A_N(i) A_N(j) - B_N(i,j)^2 \vert \frac{\eta^{\star2}}{N} + \frac{1}{2}\vert \frac{1}{2}(A_N(i)+A_N(j))\frac{\eta^\star}{\sqrt{N}} \\
  & +\frac{1}{2}(A_N(i) A_N(j) - B_N(i,j)^2) \frac{\eta^{\star2}}{N}\vert^2 \frac{3}{4}\frac{1}{\vert 1+\tilde{\beta}\vert^\frac{5}{2}},
  \end{align*}
  
  with $\vert\tilde{\beta} \vert \leq \vert \frac{1}{2}(A_N(i)+A_N(j))\frac{\eta^\star}{\sqrt{N}}+ \frac{1}{2}(A_N(i) A_N(j) - B_N(i,j)^2) \frac{\eta^{\star2}}{N} \vert \leq \epsilon_N \eta^\star +  \epsilon_N^2 \eta^{\star 2}$.
  
  $$ $$
  The remainders $\gamma_N$,$\tilde{\gamma_N}$ and $\tilde{\tilde{\gamma_N}}$ are only products of $\alpha_N$, $A_N(i)$, $A_N(j)$ and $B_N(i,j)$.

  $\vert \gamma_N \vert \leq \vert A_N(j)(A_N(i)+A_N(j))\frac{\eta^{\star2}}{N} \vert + \vert \alpha_N (1+ A_N(j)\frac{\eta^\star}{\sqrt{N}}) \vert$,

$\vert \tilde{\gamma_N} \vert \leq \vert A_N(i)(A_N(i)+A_N(j))\frac{\eta^{\star2}}{N} \vert +\vert \alpha_N (1+ A_N(i)\frac{\eta^\star}{\sqrt{N}})\vert$ and

$ \vert \tilde{\tilde{\gamma_N}}\vert \leq \vert \frac{\eta^\star}{\sqrt{N}}B_N(i,j) \vert \left( \vert (A_N(i)+A_N(j))\frac{\eta^\star}{\sqrt{N}} \vert + \vert \alpha_N \vert \right)$

 $\mu_N=\left(1 -\frac{1}{2} (A_N(i)+A_N(j)) \frac{\eta^\star}{\sqrt{N}} \right)\int_{t}^{\infty} \int_{t}^{\infty} \phi(x)\phi(y) \nu_N(x,y)dx dy$,
 
 with
 
 \begin{align*}
 \nu_N(x,y)& =-\frac{x^2}{2} \gamma_N -\frac{y^2}{2} \tilde{\gamma_N} +xy \tilde{\tilde{\gamma_N}}  \\
 & + \frac{1}{2} \left(\frac{x^2}{2}\frac{\eta^\star}{\sqrt{N}} A_N(i) +\frac{y^2}{2}\frac{\eta^\star}{\sqrt{N}} A_N(j) + xy \frac{\eta^\star}{\sqrt{N}} B_N(i,j)
  - \frac{x^2}{2} \gamma_N -\frac{y^2}{2} \tilde{\gamma_N} +xy \tilde{\tilde{\gamma_N}} \right)^2 \exp{\tilde{u}}
 \end{align*} 

$$ $$
 
 The integral of the first terms of $\nu_N(x,y)$ is 
 
 $\int_{t}^{\infty} \int_{t}^{\infty} \phi(x)\phi(y)\left(-\frac{x^2}{2} \gamma_N -\frac{y^2}{2} \tilde{\gamma_N} +xy \tilde{\tilde{\gamma_N}}\right) dxdy = -\frac{1}{2} K (t\phi(t)+K)(\gamma_N + \tilde{\gamma_N}) + \phi(t)^2 \tilde{\tilde{\gamma_N}}$.
\\
Moreover,
\\
$\exp{\tilde{u}} \leq \max(\exp\left\{\frac{x^2}{2}\frac{\eta^\star}{\sqrt{N}} A_N(i) +\frac{y^2}{2}\frac{\eta^\star}{\sqrt{N}} A_N(j) + xy \frac{\eta^\star}{\sqrt{N}} B_N(i,j) + \frac{x^2}{2} \gamma_N -\frac{y^2}{2} \tilde{\gamma_N} -xy \tilde{\tilde{\gamma_N}} \right\},1)$.
\\
There are two possibilities, either $$\max(\exp\left\{\frac{x^2}{2}\frac{\eta^\star}{\sqrt{N}} A_N(i) +\frac{y^2}{2}\frac{\eta^\star}{\sqrt{N}} A_N(j) + xy \frac{\eta^\star}{\sqrt{N}} B_N(i,j) - \frac{x^2}{2} \gamma_N -\frac{y^2}{2} \tilde{\gamma_N} +xy \tilde{\tilde{\gamma_N}} \right\},1)=1,$$
 \begin{align*}
\textrm{or }\max & (\exp\left\{\frac{x^2}{2}\frac{\eta^\star}{\sqrt{N}} A_N(i) +\frac{y^2}{2}\frac{\eta^\star}{\sqrt{N}} A_N(j) + xy \frac{\eta^\star}{\sqrt{N}} B_N(i,j) - \frac{x^2}{2} \gamma_N -\frac{y^2}{2} \tilde{\gamma_N} +xy \tilde{\tilde{\gamma_N}} \right\},1) \\
& =\exp\left\{\frac{x^2}{2}\frac{\eta^\star}{\sqrt{N}} A_N(i) +\frac{y^2}{2}\frac{\eta^\star}{\sqrt{N}} A_N(j) + xy \frac{\eta^\star}{\sqrt{N}} B_N(i,j) - \frac{x^2}{2} \gamma_N -\frac{y^2}{2} \tilde{\gamma_N} +xy \tilde{\tilde{\gamma_N}} \right\}.
\end{align*}

If $\max(\exp\left\{\frac{x^2}{2}\frac{\eta^\star}{\sqrt{N}} A_N(i) +\frac{y^2}{2}\frac{\eta^\star}{\sqrt{N}} A_N(j) + xy \frac{\eta^\star}{\sqrt{N}} B_N(i,j) - \frac{x^2}{2} \gamma_N -\frac{y^2}{2} \tilde{\gamma_N} +xy \tilde{\tilde{\gamma_N}} \right\},1)=1$,

\begin{align}
\label{eq:reste_exp}
 \int_{t}^{\infty} \int_{t}^{\infty} \phi(x)\phi(y) &\left(\frac{x^2}{2}\frac{\eta^\star}{\sqrt{N}} A_N(i) +\frac{y^2}{2}\frac{\eta^\star}{\sqrt{N}} A_N(j) + xy \frac{\eta^\star}{\sqrt{N}} B_N(i,j) - \frac{x^2}{2} \gamma_N -\frac{y^2}{2} \tilde{\gamma_N} +xy \tilde{\tilde{\gamma_N}} \right)^2 dx dy \nonumber \\
& =\frac{1}{N} J,
\end{align}
\\
where
\\
 $J=\int_{t}^{\infty} \int_{t}^{\infty} \phi(x)\phi(y) \left(\frac{x^2}{2}\eta^\star A_N(i) +\frac{y^2}{2}\eta^\star A_N(j) + xy \eta^\star B_N(i,j) - \frac{x^2}{2} \frac{\gamma_N}{\sqrt{N}} -\frac{y^2}{2}\frac{\tilde{\gamma_N}}{\sqrt{N}} +xy \frac{\tilde{\tilde{\gamma_N}}}{\sqrt{N}} \right)^2 dx dy $ is finite.
\\
 Otherwise,
 
 \begin{align*}
 \exp (\tilde{u}) \leq \exp \left\{ \frac{x^2}{2} (\epsilon_N \eta^\star +P_1(\epsilon_N)) + \frac{y^2}{2} (\epsilon_N \eta^\star+P_2(\epsilon_N)) + xy (\epsilon_N \eta^\star+P_3(\epsilon_N)) \right\}
 \end{align*}
\\ 
 where $P_1$, $P_2$, $P_3$ are polynomial functions. This expression comes from upper bounding the terms $A_N(i)/N$, $A_N(j)/N$ and $B_N(i,j)/N$ by $\epsilon_N$ in $\gamma_N$, $\tilde{\gamma}_N$ and $\tilde{\tilde{\gamma}}_N$.
 \\
 There exists $N_0$, such that for all $N \geq N_0$,  $\epsilon_N \eta^\star +P_1(\epsilon_N) \leq \frac{1}{4}$, $\epsilon_N \eta^\star +P_2(\epsilon_N) \leq \frac{1}{4}$ and $\epsilon_N \eta^\star +P_3(\epsilon_N) \leq \frac{1}{4}$.
 
 Then $\exp (\tilde{u}) \leq \exp \left\{\frac{x^2}{4} +\frac{y^2}{4} +\frac{xy}{4} \right\} \leq \exp \left\{\frac{3x^2}{8} +\frac{3y^2}{8}  \right\} $
 
 Then similarly to the expression \ref{eq:reste_exp}, 
 
 \begin{align*}
 \int_{t}^{\infty} \int_{t}^{\infty} \phi(x)\phi(y) & \left(\frac{x^2}{2}\frac{\eta^\star}{\sqrt{N}} A_N(i) +\frac{y^2}{2}\frac{\eta^\star}{\sqrt{N}} A_N(j) + xy \frac{\eta^\star}{\sqrt{N}} B_N(i,j) + \frac{x^2}{2} \gamma_N -\frac{y^2}{2} \tilde{\gamma_N} -xy \tilde{\tilde{\gamma_N}} \right)^2 \exp(\tilde{u})dx dy \\
& \leq \frac{1}{2\pi}\int_{t}^{\infty} \int_{t}^{\infty} \exp(-\frac{x^2}{8})\exp(-\frac{y^2}{8}) \left(\frac{x^2}{2}\frac{\eta^\star}{\sqrt{N}} A_N(i) \right. \\
& \left. +\frac{y^2}{2}\frac{\eta^\star}{\sqrt{N}} A_N(j) + xy \frac{\eta^\star}{\sqrt{N}} B_N(i,j) + \frac{x^2}{2} \gamma_N -\frac{y^2}{2} \tilde{\gamma_N} -xy \tilde{\tilde{\gamma_N}} \right)^2 dxdy \\
& \leq \frac{1}{N} J'
\end{align*}
where $J'$ is finite.

 Similarly to the computations made for $\alpha_N$, $\beta_N$, $\gamma_N$, $\nu_N$, all the remainder terms can be upper bounded by products of $A_N(i)/\sqrt{N}$, $A_N(j)/\sqrt{N}$ and  $B_N(i,j)/\sqrt{N}$, which proves \eqref{eq:proof_RN}.

%
%
%
%
%
%

\subsubsection{Proof of Lemma \ref{lemma4}}
\label{sec:proof_4}
In this section, all the expectations that we consider are conditionally to the presence of the observed individuals in the study, for instance $\left\{\epsilon_i=\epsilon_j=1\right\}$ or $\left\{\epsilon_{i_1}=\epsilon_{i_2}=\epsilon_{i_3}=1\right\}$. However, for the sake of simplicity, we will not always make explicit such conditioning.
\\
 Let us show that $$\Var(\frac{1}{n}\underset{i \neq j}{\sum}( \W_i \W_j - \mathbb{E}[\W_i \W_j \vert \Z])\G_N(i,j)\mathbbm{1}_{E_N}) \to 0,$$
 \\
 that is 
 
\begin{equation}
\frac{1}{n^2}\underset{i_3 \neq i_4}{\underset{i_1 \neq i_2}{\sum}} \mathbb{E} \left[( \mathbb{E}[\W_{i_1} \W_{i_2} \W_{i_3} \W_{i_4}\vert \Z] - \mathbb{E}[\W_{i_1} \W_{i_2}\vert \Z]\mathbb{E}[\W_{i_3} \W_{i_4}\vert \Z])\G_N(i_1,i_2) \G_N(i_3,i_4)\mathbbm{1}_{E_N}\right] \to 0
\label{eq:var_4termes}
\end{equation}
\\
For this purpose, we will separate three cases depending on the cardinal of the set $\{i_1,i_2,i_3,i_4\}$ in the sum of Equation \eqref{eq:var_4termes}.

\bigskip 

-If card($\{i_1,i_2,i_3,i_4\}$)=2, the corresponding terms in \eqref{eq:var_4termes} are equal to

\begin{align*}
\frac{1}{n^2}\underset{i \neq j}{\sum}\mathbb{E} \left[ \mathbb{E}[( \W_i^2 \W_j^2 \vert \Z] - \mathbb{E}[\W_i \W_j \vert \Z]^2)\G_N(i,j)^2 \mathbbm{1}_{E_N} \right] \leq \frac{1}{n^2}\underset{i \neq j}{\sum} \mathbb{E}[(\alpha +   \rho_N(i,j) ) \G_N(i,j)^2\mathbbm{1}_{E_N}]
\end{align*}
\\
 where $\alpha$ is a positive constant and $\rho_N(i,j)$ can be upper bounded by a finite product of $\G_N(i,j)$, $\G_N(i,i)-1$ and $\G_N(j,j)-1$, according to proof of Lemma \ref{lemma3}.
This result is obtained by using a similar decomposition of $\mathbb{E}[ \W_i^2 \W_j^2 \vert \Z]$ than the one that we explicited for $\mathbb{E}[ \W_i \W_j \vert \Z]$.
\\
Since $\mathbb{E}[\G_N(i,j)^2 \mathbbm{1}_{E_N}]\leq \epsilon_N^2$ and all terms of $ \rho_N(i,j)$ are upper bounded by a finite sum of $\epsilon_N^k$, with $k$ greater than 1, which all tend to $0$, it is clear that

 $$\frac{1}{n^2}\underset{i \neq j}{\sum}\mathbb{E} \left[ \mathbb{E}[( \W_i^2 \W_j^2 \vert \Z] - \mathbb{E}[\W_i \W_j \vert Z]^2 )\G_N(i,j)^2 \mathbbm{1}_{E_N}\right] \to 0.$$
\\
- If card($\{i_1,i_2,i_3,i_4\}$)=3, the corresponding terms in \eqref{eq:var_4termes} are equal to
\\
\begin{equation}
\frac{1}{n^2}\underset{i_1 \neq i_2 \neq i_3}{\sum} \mathbb{E} \left[( \mathbb{E}[\W_{i_1}^2 \W_{i_2} \W_{i_3} \vert \Z] - \mathbb{E}[\W_{i_1} \W_{i_2}\vert Z]\mathbb{E}[\W_{i_1} \W_{i_3}\vert Z])\G_N(i_1,i_2) \G_N(i_1,i_3) \mathbbm{1}_{E_N}\right].
\label{eq:var_4termes_card3}
\end{equation}
\\
Since the sum of Equation \eqref{eq:var_4termes_card3} has $n(n-1)(n-2)$ terms, we have the refine the upper bound that we used in the case where the cardinal of $\{i_1,i_2,i_3,i_4\}$ was equal to $2$. Indeed, we will use the following proposition:

$$ $$
\begin{prop}
$  \mathbb{E}[\W_{i_1}^2 \W_{i_2} \W_{i_3} \vert \Z]$ has no term of order less than $1/\sqrt{N}$, that is no constant term.
\label{prop:no_const_term}
\end{prop}

$$ $$
\\
Let us explain why Proposition \ref{prop:no_const_term} is enough to prove 
\\
 \begin{equation}
\frac{1}{n^2}\underset{i_1 \neq i_2 \neq i_3}{\sum} \mathbb{E} \left[( \mathbb{E}[\W_{i_1}^2 \W_{i_2} \W_{i_3} \vert \Z] - \mathbb{E}[\W_{i_1} \W_{i_2}\vert \Z]\mathbb{E}[\W_{i_1} \W_{i_3}\vert \Z])\G_N(i_1,i_2) \G_N(i_1,i_3) \mathbbm{1}_{E_N}\right] \to 0.
\label{eq:pr_card3}
 \end{equation}
 \\
Let us first recall that, according to Lemma \ref{lemma3},
\\
\begin{align*}
\mathbb{E}[\W_{i_1} \W_{i_2}\vert \Z]\mathbb{E}[\W_{i_1} \W_{i_3}\vert \Z]&= c^2 \eta^{\star 2} \G_N(i_1,i_2)\G_N(i_1,i_3) +c\eta^\star \G_N(i_1,i_3) R_N(i_1,i_2) \\
&+c\eta^\star \G_N(i_1,i_2) R_N(i_1,i_3) +R_N(i_1,i_2) R_N(i_1,i_3),
\end{align*}
\\
where, if $E_N$ holds,  all these terms are upper bounded by a finite sum of terms of the form $\epsilon_N^k$, with $k\geq 2$.
\\
Then, 
$$\mathbb{E}\left[\mathbb{E}[\W_{i_1} \W_{i_2}\vert \Z]\mathbb{E}[\W_{i_1} \W_{i_3}\vert \Z]\G_N(i_1,i_2) \G_N(i_1,i_3)\mathbbm{1}_{E_N} \right] $$
can be upper bounded by a finite sum of  terms of the form $\epsilon_N^k$, with $k\geq 4$.
\\
Since 
$$ \frac{N(N-1)(N-2) \epsilon_N^4}{n^2} \to 0,$$
\\
it shows that 
\\
$$\frac{1}{n^2}\underset{i_1 \neq i_2 \neq i_3}{\sum}\mathbb{E}\left[\mathbb{E}[\W_{i_1} \W_{i_2}\vert \Z]\mathbb{E}[\W_{i_1} \W_{i_3}\vert \Z])\mathbbm{1}_{E_N} \right] \to 0.$$

Similarly, according to Proposition \ref{prop:no_const_term}, 
 each term of $\mathbb{E}\left[ \mathbb{E}[\W_{i_1}^2 \W_{i_2} \W_{i_3} \vert \Z] \G_N(i_1,i_2) \G_N(i_1,i_3) \mathbbm{1}_{E_N}\right]$ can be upper bounded by a finite sum of $\epsilon_N^k$, with $k\geq 3$. 
 \\
 Since 
 
 $$\frac{n(n-1)(n-2) \epsilon_N^3}{n^2}= O \left( \frac{1}{N^{1/2 - 3\gamma}} \right) \to 0,$$
%
\\
it achieves the proof of \eqref{eq:var_4termes_card3}.
%
%
%
%
%
\\
- If card($\{i_1,i_2,i_3,i_4\}$)=4, let us first observe that 

$$ \frac{N(N-1)(N-2)(N-3) \epsilon_N^5}{n^2} \to 0,$$
\\
which means that we shall only focus on the approximation of 
\\
$  \mathbb{E}[\W_{i_1} \W_{i_2} \W_{i_3} \W_{i_4}\vert \Z] - \mathbb{E}[\W_{i_1} \W_{i_2}\vert \Z]\mathbb{E}[\W_{i_3} \W_{i_4}\vert \Z]$ of order $1/N$.
%
%
%
\\
Let us recall that 
\\
$$\mathbb{E}\left[\W_{i_1} \W_{i_2}\vert \Z] \mathbb{E}[\W_{i_3} \W_{i_4}\vert \Z\right]  = c^2 \eta^{\star2} \G_N(i_1,i_2)\G_N(i_3,i_4) + R_N(i_1,i_2,i_3,i_4), $$
 where $$R_N(i_1,i_2,i_3,i_4)=c\eta^\star \G_N(i_1,i_2)R_N(i_3,i_4) + c\eta^\star \G_N(i_3,i_4) R_N(i_1,i_2)+ R_N(i_1,i_2)R_N(i_3,i_4)$$ is a remainder, each term of which is upper bounded by a finite sum of terms of the form $\epsilon_N^k$, with $k\geq 2$. In particular, it implies that $$ \mathbb{E}[ \frac{N(N-1)(N-2)(N-3)}{n^2} R_N(i_1,i_2,i_3,i_4) \G_N(i_1,i_2)\G_N(i_3,i_4)] \to 0.$$
\\
Thus, we need to prove that 
 \begin{equation}
\frac{1}{n^2}\underset{i_1 \neq i_2 \neq i_3}{\sum} \mathbb{E} \left[( \mathbb{E}[\W_{i_1} \W_{i_2} \W_{i_3} \W_{i_4} \vert \Z] - c^2 \eta^{\star2} \G_N(i_1,i_2)\G_N(i_3,i_4))\G_N(i_1,i_2) \G_N(i_3,i_4) \mathbbm{1}_{E_N}\right] \to 0,
\label{eq:pr_card4_eq}
 \end{equation}

$$ $$
To do so, we shall prove first the following proposition:
$$ $$
\begin{prop}
The terms of order less than or equal to $1/\sqrt{N}$ in $\mathbb{E}[\W_{i_1} \W_{i_2} \W_{i_3} \W_{i_4}\vert \Z]$ are null.
\label{prop1_card4}
\end{prop}
$$ $$
The term of order exactly $1/N$ in $\mathbb{E}[\W_{i_1} \W_{i_2} \W_{i_3} \W_{i_4}\vert \Z]$ contains all combinations of products of two terms between $\G_N(i_1,i_2)$, $\G_N(i_1,i_3)$, $\G_N(i_1,i_4)$, $\G_N(i_2,i_3)$, $\G_N(i_2,i_4)$, $\G_N(i_3,i_4)$,  $\G_N(i_1,i_1)-1$, $\G_N(i_2,i_2)-1$, $\G_N(i_3,i_3)-1$ and $\G_N(i_4,i_4)-1$.
\\
We will demonstrate the propositions:
\\
\begin{prop}
The term in $\G_N(i_1,i_2)\G_N(i_3,i_4)$ of $\mathbb{E}[\W_{i_1} \W_{i_2} \W_{i_3} \W_{i_4}\vert \Z]$ is equal to 
$c^2 \eta^{\star2} \G_N(i_1,i_2)\G_N(i_3,i_4)$.
\label{prop3_card4}
\end{prop}
\bigskip
\begin{prop}
For all terms $T_N(i_1,i_2,i_3,i_4)$ of order $1/N$ in $\mathbb{E}[\W_{i_1} \W_{i_2} \W_{i_3} \W_{i_4}\vert \Z]$, 
\\
$$ \frac{1}{n^2}\mathbb{E}[T_N(i_1,i_2,i_3,i_4) \G_N(i_1,i_2)\G_N(i_3,i_4) ]\to 0,$$
 except for the term in $\G_N(i_1,i_2)\G_N(i_3,i_4)$.
 \label{prop2_card4}
\end{prop}
Propositions \ref{prop1_card4}, \ref{prop3_card4} and \ref{prop2_card4} prove \eqref{eq:pr_card4_eq}.
\\
Let us prove now  Propositions \ref{prop:no_const_term}, \ref{prop1_card4}, \ref{prop2_card4} and \ref{prop3_card4}.

$$ $$
\\
If card($\{i_1,i_2,i_3,i_4\}$)=3,
conditionally to $\left\{\epsilon_{i_1}= \epsilon_{i_2}=\epsilon_{i_3}=1 \right\}$,
$\W_{i_1}^2 \W_{i_2} \W_{i_3}$ can take several values:

$\bullet$ $\frac{(1-P)^2}{P^2}$ if $\Y_{i_1}=\Y_{i_2}=\Y_{i_3}=1$.

$\bullet$ $\frac{-(1-P)}{P}$ if $\Y_{i_1}=1$ and $\Y_{i_2}\neq \Y_{i_3}$.

$\bullet$ $1$ if $\Y_{i_1}=1$ and $\Y_{i_2}=\Y_{i_3}=0$ or $\Y_{i_1}=0$ and $\Y_{i_2}=\Y_{i_3}=1$.

$\bullet$ $\frac{-P}{1-P}$ if $\Y_{i_1}=0$ and $\Y_{i_2}\neq \Y_{i_3}$.

$\bullet$ $\frac{P^2}{(1-P)^2}$ if $\Y_{i_1}=\Y_{i_2}=\Y_{i_3}=0$.
\\
Since each case has a probability $1$ and each control a probability $K(1-P)/P(1-K)$ to be in the study (these probabilities are given in Equation \eqref{eq:full_asc} and \eqref{eq:p_control}),

\begin{align*}
&\mathbb{E}[\W_{i_1}^2 \W_{i_2} \W_{i_3} \vert \Z, \epsilon_{i_1}=\epsilon_{i_2}=\epsilon_{i_3}=1]=\frac{1}{\mathbb{P}(\epsilon_{i_1}=\epsilon_{i_2}=\epsilon_{i_3}=1)} \\ & \times \left\{\frac{(1-P)^2}{P^2} \mathbb{P}(\Y_{i_1}=\Y_{i_2}=\Y_{i_3}=1 \vert \Z)  -  \frac{1-P}{P}\left(\frac{K(1-P)}{P(1-K)} \right)\mathbb{P}(\Y_{i_1}=1, \Y_{i_2}\neq \Y_{i_3} \vert \Z) \right. \\
& \left. + \left(\frac{K(1-P)}{P(1-K)} \right) \mathbb{P}(\Y_{i_1}=0,\Y_{i_2}=\Y_{i_3}=1 \vert \Z)+ \left(\frac{K(1-P)}{P(1-K)} \right)^2 \mathbb{P}(\Y_{i_1}=1,\Y_{i_2}=\Y_{i_3}=0 \vert \Z) \right. \\
& \left. - \frac{P}{1-P}\left(\frac{K(1-P)}{P(1-K)} \right)^2 \mathbb{P}(\Y_{i_1}=0, \Y_{i_2}\neq \Y_{i_3} \vert \Z) + \frac{(1-P)^2}{P^2}\left(\frac{K(1-P)}{P(1-K)} \right)^3\mathbb{P}(\Y_{i_1}=\Y_{i_2}= \Y_{i_3}=0\vert \Z)  \right\}
\end{align*}
\\
The development of order 0 of $\mathbb{P}(\Y_{i_1}=\Y_{i_2}=\Y_{i_3}=1 \vert \Z)$ is 

$$\frac{1}{(2\pi)^{\frac{3}{2}}}\int_t^{+\infty} \int_t^{+\infty}\int_t^{+\infty} \phi(x) \phi(y) \phi(z) dx dy dz = K^3 + O_p\left(\frac{1}{\sqrt N}\right).$$
\\
Similarly, 

$$\mathbb{P}(\Y_{i_1}=1, \Y_{i_2}\neq \Y_{i_3} \vert \Z) = 2K^2(1-K)+O_p\left(\frac{1}{\sqrt N}\right)$$
$$\mathbb{P}(\Y_{i_1}=0, \Y_{i_2}= \Y_{i_3}=1 \vert \Z) = K^2(1-K)+O_p\left(\frac{1}{\sqrt N}\right)$$
$$\mathbb{P}(\Y_{i_1}=1, \Y_{i_2}= \Y_{i_3}=0 \vert \Z) = K(1-K)^2+O_p\left(\frac{1}{\sqrt N}\right)$$
$$\mathbb{P}(\Y_{i_1}=0, \Y_{i_2}\neq \Y_{i_3} \vert \Z) = 2K(1-K)^2+O_p\left(\frac{1}{\sqrt N}\right)$$
$$\mathbb{P}(\Y_{i_1}=\Y_{i_2}=\Y_{i_3}=0 \vert \Z) =(1-K)^3 +O_p\left(\frac{1}{\sqrt N}\right) $$
\\
Replacing all these expressions in $\mathbb{E}[\W_{i_1}^2 \W_{i_2} \W_{i_3} \vert \Z, \epsilon_{i_1}=\epsilon_{i_2}=\epsilon_{i_3}=1]$ gives us that the approximation of order 0 is null, which achieves the proof of Proposition \ref{prop:no_const_term}.
%
%
%
%
%
%
%
%
%
\\
Let us prove now Proposition \ref{prop1_card4}.
\\
If card($\{i_1,i_2,i_3,i_4\}$)=4, let us compute the approximation of order $1/\sqrt{N}$ of $\mathbb{E}[\W_{i_1} \W_{i_2} \W_{i_3} \W_{i_4}\vert \Z]$.
\\
Conditionally to $\left\{\epsilon_{i_1}=\epsilon_{i_2}=\epsilon_{i_3}=\epsilon_{i_4}=1\right\}$, 
$\W_{i_1} \W_{i_2} \W_{i_3} \W_{i_4}$ can take values:

$\bullet$ $\frac{(1-P)^2}{P^2}$ if all individuals are cases, that is $\Y_{i_1}=\Y_{i_2}=\Y_{i_3}=\Y_{i_4}=1$.

$\bullet$  $\frac{-(1-P)}{P}$ if one individual is a control and the three others are cases.

$\bullet$ $1$ if two individuals are controls and two are cases.

$\bullet$  $\frac{-P}{1-P}$ if one individual is a case and the three others are controls.

$\bullet$ $\frac{P^2}{(1-P)^2}$ if all individuals are controls.
\\
The matrix of variance covariance of $(\l_{i_1},\l_{i_2},\l_{i_3},\l_{i_4})$ is
\\
$$ \Sigma = \begin{pmatrix}
1+\eta^\star (\G_N(i_1,i_1)-1) &\G_N(i_1,i_2) & \G_N(i_1,i_3)&\G_N(i_1,i_4) \\
\G_N(i_1,i_2) & 1+\eta^\star (\G_N(i_2,i_2)-1) & \G_N(i_2,i_3) & \G_N(i_2,i_4) \\
\G_N(i_1,i_3) &\G_N(i_2,i_3)&1+\eta^\star (\G_N(i_3,i_3)-1)& \G_N(i_3,i_4)\\
\G_N(i_1,i_4)&\G_N(i_2,i_4)&\G_N(i_3,i_4) & 1+\eta^\star (\G_N(i_4,i_4)-1)
\end{pmatrix}.
$$
\\
For the sake of clarity, let us denote $A_1=  \frac{1}{\sqrt N} \overset{N}{\underset{k=1}{\sum}}( \Z_{i_1,k}^2-1) = \sqrt N(\G_N(i_1,i_1)-1)$, and similarly we define $A_2$, $A_3$ and $A_4$.
\\
Let us also denote $C_{1,2}=  \sqrt N \G_N(i_1,i_2)$ and similarly, $C_{1,3}, \dots, C_{3,4}$.
\\
Then, let us rewrite $\Sigma$ as: 

$$ \Sigma = \begin{pmatrix}
1+\frac{\eta^\star}{\sqrt N} A_1&\frac{C_{1,2}}{\sqrt n} & \frac{C_{1,3}}{\sqrt N}&\frac{C_{1,4}}{\sqrt N} \\
\frac{C_{1,2}}{\sqrt N} &1+\frac{\eta^\star}{\sqrt N} A_2& \frac{C_{2,3}}{\sqrt N} & \frac{C_{2,4}}{\sqrt N} \\
\frac{C_{1,3}}{\sqrt N} &\frac{C_{2,3}}{\sqrt N}&1+\frac{\eta^\star}{\sqrt N} A_3& \frac{C_{3,4}}{\sqrt N}\\
\frac{C_{1,4}}{\sqrt N}&\frac{C_{2,4}}{\sqrt N}&\frac{C_{3,4}}{\sqrt N} & 1+\frac{\eta^\star}{\sqrt N} A_4
\end{pmatrix}.
$$
\\
The approximation of order $1/\sqrt n$ of its inverse matrix is given by 

\begin{align*}
 \Sigma ^{-1}& \simeq \vert \Sigma \vert ^{-1} \\
 & \times  \begin{pmatrix}
1+\frac{\eta^\star}{\sqrt N} (A_2+A_3+A_4) &-\frac{C_{1,2}}{\sqrt N} &- \frac{C_{1,3}}{\sqrt N}&-\frac{C_{1,4}}{\sqrt N} \\
-\frac{C_{1,2}}{\sqrt N} &1+\frac{\eta^\star}{\sqrt N} (A_1+A_3+A_4) & -\frac{C_{2,3}}{\sqrt N} &- \frac{C_{2,4}}{\sqrt N} \\
-\frac{C_{1,3}}{\sqrt N} &-\frac{C_{2,3}}{\sqrt N}&1+\frac{\eta^\star}{\sqrt N} (A_1+A_2+A_4) &- \frac{C_{3,4}}{\sqrt N}\\
-\frac{C_{1,4}}{\sqrt N}&-\frac{C_{2,4}}{\sqrt N}&-\frac{C_{3,4}}{\sqrt N} & 1+\frac{\eta^\star}{\sqrt n}(A_1+A_2+A_2) 
\end{pmatrix}
\end{align*}
\\
where $\vert \Sigma \vert ^{-1} = 1 - \frac{\eta^\star}{\sqrt N}(A_1+A_2+A_3+A_4) + O_p\left(\frac{1}{N}\right)$.
\\
Let us compute

$$\mathbb{P}(\Y_{i_1}=\Y_{i_2}=\Y_{i_3}=\Y_{i_4}=1 \vert \Z)=  \frac{1} {\vert \Sigma \vert^\frac{1}{2}} \int_t^{+\infty}\int_t^{+\infty}\int_t^{+\infty}\int_t^{+\infty} f(w,x,y,z) dw dx dy dz,$$
\\
where 

\begin{align*}
\hspace{-10mm} f(w,x,y,z)& = \frac{1}{(2\pi)^2} \exp \left\{ - \frac{x^2}{2 \vert \Sigma \vert} (1+\frac{\eta^\star}{\sqrt n} (A_2+A_3+A_4)) - ...   - \frac{z^2}{2 \vert \Sigma \vert} (1+\frac{\eta^\star}{\sqrt N} (A_1+A_2+A_3)) \right. \\
 & \left.+ \frac{wx}{ \vert \Sigma \vert} \frac{\eta^\star}{\sqrt n} C_{1,2} + \frac{wy}{ \vert \Sigma \vert} \frac{\eta^\star}{\sqrt n}C_{1,3}+ ... + \frac{yz}{ \vert \Sigma \vert} \frac{\eta^\star}{\sqrt N} C_{3,4} + O_p\left(\frac{1}{N}\right) \right\} \\
 &= \frac{1}{(2\pi)^2} \exp \left\{ - \frac{x^2}{2} (1-\frac{\eta^\star}{\sqrt N}(A_1+A_2+A_3+A_4)) (1+\frac{\eta^\star}{\sqrt N} (A_2+A_3+A_4)) - ...     \right.\\
 & - \frac{z^2}{2}(1-\frac{\eta^\star}{\sqrt N}(A_1+A_2+A_3+A_4)) (1+\frac{\eta^\star}{\sqrt n} (A_1+A_2+A_3))  \\
 & \left.+ wx \frac{\eta^\star}{\sqrt N} (1-\frac{\eta^\star}{\sqrt N}(A_1+A_2+A_3+A_4)) C_{1,2} + ... + yz(1-\frac{\eta^\star}{\sqrt N}(A_1+A_2+A_3+A_4)) \frac{\eta^\star}{\sqrt N} C_{3,4} \right\} \\
 &= \phi(w) \phi(x) \phi(y) \phi(z) \exp \left\{ \frac{x^2}{2} \frac{\eta^\star}{\sqrt N} A_1 + ... +\frac{z^2}{2} \frac{\eta^\star}{\sqrt N} A_4 - wx \frac{\eta^\star}{\sqrt N} C_{1,2} - ... - yz \frac{\eta^\star}{\sqrt N} C_{3,4} + O_p\left(\frac{1}{N}\right) \right\} \\
 &=  \phi(w) \phi(x) \phi(y) \phi(z) \left[ 1+\frac{x^2}{2} \frac{\eta^\star}{\sqrt N} A_1 + ... +\frac{z^2}{2} \frac{\eta^\star}{\sqrt N} A_4 - wx \frac{\eta^\star}{\sqrt N} C_{1,2} - ... - yz \frac{\eta^\star}{\sqrt N} C_{3,4} + O_p\left(\frac{1}{N}\right)\right]
\end{align*}
\\
Finally,
\begin{align*}
\bullet \textrm{ }  \mathbb{P}(\Y_{i_1}=\Y_{i_2}=\Y_{i_3}=\Y_{i_4}=1  \vert \Z)&= \frac{1}{ \vert \Sigma \vert^\frac{1}{2}} \left[ K^4 + \frac{K^3}{2}(t\phi(t)+K)  \frac{\eta^\star}{\sqrt n} (A_1+A_2+A_3+A_4) \right.  \\
& \left.+ K^2\phi(t)^2  \frac{\eta^\star}{\sqrt n} (C_{1,2} + ... + C_{3,4}) \right]
\end{align*}
\\
Similarly, we compute 
\\
\begin{align*}
\bullet \textrm{ }   \mathbb{P}(\textrm{``1 control, 3 cases"}\vert \Z)&=\frac{1}{ \vert \Sigma \vert^\frac{1}{2}} \left[4K^3(1-K)  \right. \\
& + \frac{K^2}{2} \left((3-4K) t\phi(t) + 4K(1-K) \right)\frac{\eta^\star}{\sqrt N} (A_1 + A_2 + A_3+A_4) \\
& \left. + 2 \phi(t)^2 K(1-2K) \frac{\eta^\star}{\sqrt N} (C_{1,2} + ... + C_{3,4}) \right] 
\end{align*}
\\
\begin{align*}
\bullet \textrm{ }   \mathbb{P}(\textrm{``2 controls, 2 cases"}\vert \Z)&=\frac{1}{ \vert \Sigma \vert^\frac{1}{2}} \left[6K^2(1-K)^2  \right. \\
& + \frac{3 K(1-K)}{2} \left((1-2K) t\phi(t) + 2K(1-K) \right)\frac{\eta^\star}{\sqrt N} (A_1 + A_2 + A_3+A_4) \\
& \left. + \phi(t)^2 (6 K^2 -6K +1) \frac{\eta^\star}{\sqrt N} (C_{1,2} + ... + C_{3,4}) \right] 
\end{align*}
\\
\begin{align*}
\bullet \textrm{ }   \mathbb{P}(\textrm{``3 controls, 1 case"}\vert \Z)&=\frac{1}{ \vert \Sigma \vert^\frac{1}{2}} \left[4K(1-K)^3  \right. \\
& + \frac{(1-K)^2}{2} \left((1-4K) t\phi(t) + 4K(1-K) \right)\frac{\eta^\star}{\sqrt n} (A_1 + A_2 + A_3+A_4) \\
& \left. - 2 \phi(t)^2 (1-K)(1-2K) \frac{\eta^\star}{\sqrt N} (C_{1,2} + ... + C_{3,4}) \right] 
\end{align*}
\\
\begin{align*}
\bullet \textrm{ }  \mathbb{P}(\Y_{i_1}=\Y_{i_2}=\Y_{i_3}=\Y_{i_4}=0  \vert \Z)&= \frac{1}{ \vert \Sigma \vert^\frac{1}{2}} \left[ (1-K)^4 \right.\\
&+ \frac{(1-K)^3}{2}(-t\phi(t)+1-K)  \frac{\eta^\star}{\sqrt n} (A_1+A_2+A_3+A_4)  \\
& \left.+ (1-K)^2\phi(t)^2  \frac{\eta^\star}{\sqrt N} (C_{1,2} + ... + C_{3,4}) \right]
\end{align*}
\\
 Regrouping all the first terms in the expression of $\mathbb{E}[\W_{i_1} \W_{i_2} \W_{i_3} \W_{i_4}\vert \Z] $ gives
\begin{align*}
\frac{1}{ \vert \Sigma \vert^\frac{1}{2}}&  \left[ \frac{(1-P)^2}{P^2} K^4 - \frac{(1-P)}{P} \left( \frac{K(1-P)}{P(1-K)} \right) 4K^3(1-K) + \left( \frac{K(1-P)}{P(1-K)} \right) ^2 6K^2(1-K)^2 \right. \\
& \left. - \frac{P}{1-P} \left( \frac{K(1-P)}{P(1-K)} \right)^3 4K(1-K)^3 + \frac{P^2}{(1-P)^2}\left( \frac{K(1-P)}{P(1-K)} \right)^4 (1-K)^4 \right] \\
& = \frac{1}{ \vert \Sigma \vert^\frac{1}{2}} \left( \frac{(1-P)^2 K^4}{P^2}\right) \left[ 1 - 4 + 6 - 4+ 1 \right] =0
\end{align*}
\\
Similarly we regroup the terms in $\frac{\eta^\star}{\sqrt N} (A_1 + A_2 + A_3+A_4)$: 

\begin{align*}
\frac{1}{  \vert \Sigma \vert^\frac{1}{2}}& \frac{\eta^\star}{\sqrt N} (A_1 + A_2 + A_3+A_4)  \left[ \frac{(1-P)^2}{P^2}  \frac{K^3}{2}(t\phi(t)+K)  \right. \\
& - \frac{(1-P)}{P} \left( \frac{K(1-P)}{P(1-K)} \right)  \frac{K^2}{2} \left((3-4K) t\phi(t) + 4K(1-K) \right) \\
& + \left( \frac{K(1-P)}{P(1-K)} \right) ^2 \frac{3 K(1-K)}{2} \left((1-2K) t\phi(t) + 2K(1-K) \right)  \\
& - \frac{P}{1-P} \left( \frac{K(1-P)}{P(1-K)} \right)^3 \frac{(1-K)^2}{2} \left((1-4K) t\phi(t) + 4K(1-K) \right) \\
& \left.+ \frac{P^2}{(1-P)^2}\left( \frac{K(1-P)}{P(1-K)} \right)^4 \frac{(1-K)^3}{2}(-t\phi(t)+1-K) \right] \\
& = \frac{1}{ \vert \Sigma \vert^\frac{1}{2}} \left( \frac{(1-P)^2 K^4}{2 P^2}\right) \left[ 1 - 4 + 6 - 4+ 1 \right]  \\
& + \frac{1}{ \vert \Sigma \vert^\frac{1}{2}} \left( \frac{(1-P)^2 K^3}{2 P^2 (1-K)}\right)  \left[1-K - 3 +4K + 3(1-2K) -1 +4K -K\right]=0
\end{align*}
\\
Finally, we regroup all the terms in $ \frac{\eta^\star}{\sqrt n} (C_{1,2} + ... + C_{3,4})$:

\begin{align*}
 \frac{1}{ \vert \Sigma \vert^\frac{1}{2}} \left( \frac{(1-P)^2 K^2}{ P^2(1-K)^2}\right)  \phi(t)^2 & \left[ (1-K)^2 - 2(1-K)(1-2K) + 6K^2-6K+1 +2K(1-2K)+K^2 \right]  \\
&=0.
\end{align*}
\\
This proves Proposition \ref{prop1_card4}.

$$ $$
Let us prove Proposition \ref{prop2_card4}.
\\
The main term of the second order approximation of $f(w,x,y,z)$ can be written as:

\begin{align}
f_2(w,x,y,z)& = \phi(w) \phi(x) \phi(y) \phi(z) \left[1 + \frac{w^2}{2} (\frac{\eta^\star}{\sqrt N} A_1- \frac{\eta^{\star 2}}{N}(A_1^2 + C_{1,2}^2 +C_{1,3}^2+C_{1,4}^2) + \dots \right. \nonumber \\
& + \frac{z^2}{2} (\frac{\eta^\star}{\sqrt N} A_4- \frac{\eta^{\star 2}}{N}(A_4^2 + C_{1,4}^2 +C_{2,4}^2+C_{3,4}^2) \nonumber\\
& + wx (C_{1,2} \frac{\eta^\star}{\sqrt N} - \frac{\eta^{\star2}}{N} [(A_1+A_2) C_{1,2} + C_{1,3}C_{2,3} + C_{1,4}C_{2,4}]) +  \dots \nonumber\\
& + yz (C_{3,4} \frac{\eta^\star}{\sqrt N} - \frac{\eta^{\star2}}{N} [(A_3+A_4) C_{3,4} + C_{1,3}C_{1,4} + C_{2,3}C_{2,4}])\nonumber \\
&+ \frac{w^4}{8} \frac{\eta^{\star2}}{N} A_1^2+ \dots+ \frac{z^4}{8} \frac{\eta^{\star2}}{N} A_4^2 + \frac{w^2x^2}{2} \frac{\eta^{\star2}}{N} (C_{1,2}^2+\frac{A_1A_2}{2})+ \dots \nonumber\\
& +  \frac{y^2z^2}{2} \frac{\eta^{\star2}}{N} (C_{3,4}^2+\frac{A_3A_4}{2}) + \frac{w^3x}{2}\frac{\eta^{\star2}}{N} A_1 C_{1,2} + … \nonumber\\
& + \frac{z^3y}{2}\frac{\eta^{\star2}}{N} A_4 C_{3,4} + w^2xy \frac{\eta^{\star 2}}{N} [ \frac{A_1 C_{2,3}}{2} + C_{1,2}C_{1,3} ] +  \dots\nonumber \\
& \left. + z^2xy \frac{\eta^{\star 2}}{N}[ \frac{A_4 C_{2,3}}{2} + C_{2,4}C_{3,4} ] + wxyz \frac{\eta^{\star 2}}{N}(C_{1,2}C_{3,4} + C_{2,3}C_{1,4} + C_{1,3}C_{2,4}) \right].
\label{eq:f2}
\end{align}
\\
In order to prove Proposition \ref{prop2_card4}, we will show that: 

\begin{equation}
\frac{1}{n^2} \underset{i_1\neq i_2 \neq i_3 \neq i_4}{\sum} \mathbb{E}[A_1^2 C_{1,2}C_{3,4}] \to 0
\label{eq:A1_2}
\end{equation}

\begin{equation}
\frac{1}{n^2} \underset{i_1\neq i_2 \neq i_3 \neq i_4}{\sum} \mathbb{E}[A_1A_2 C_{1,2}C_{3,4}] \to 0
\label{eq:A1_A2}
\end{equation}

\begin{equation}
\frac{1}{n^2} \underset{i_1\neq i_2 \neq i_3 \neq i_4}{\sum} \mathbb{E}[A_1 C_{1,2}^2C_{3,4}] \to 0
\label{eq:A1_C12}
\end{equation}

\begin{equation}
\frac{1}{n^2} \underset{i_1\neq i_2 \neq i_3 \neq i_4}{\sum} \mathbb{E}[A_1 C_{1,2} C_{13} C_{3,4}] \to 0
\label{eq:A1_C13}
\end{equation}

\begin{equation}
\frac{1}{n^2} \underset{i_1\neq i_2 \neq i_3 \neq i_4}{\sum} \mathbb{E}[C_{1,2}^2 C_{2,3} C_{3,4}] \to 0
\label{eq:C12_C23}
\end{equation}

\begin{equation}
\frac{1}{n^2} \underset{i_1\neq i_2 \neq i_3 \neq i_4}{\sum} \mathbb{E}[C_{1,2} C_{1,3} C_{2,4} C_{3,4}] \to 0
\label{eq:C13_C24}
\end{equation}

\begin{equation}
\frac{1}{n^2} \underset{i_1\neq i_2 \neq i_3 \neq i_4}{\sum} \mathbb{E}[C_{1,2}^3 C_{3,4}] \to 0
\label{eq:C12_2}
\end{equation}
\\
We will develop the proof of Equation \eqref{eq:A1_A2}.
\\
By exchangeability of the $(\Z_{i,k})_{1\leq i \leq n}$, we can write

\begin{align}
\mathbb{E}[A_1 A_2 C_{1,2} C_{3,4}] & = \underset{k,l,m,r}{\sum} \mathbb{E}[(\Z_{1,k}^2-1)(\Z_{2,l}^2-1)\Z_{1,m}\Z_{2,m}\Z_{3,r} \Z_{4,r} ] \\
&= \underset{k,l,m,r}{\sum} \mathbb{E}[\Z_{1,k}^2\Z_{2,l}^2\Z_{1,m}\Z_{2,m}\Z_{3,r} \Z_{4,r} ] - 2N \underset{k,m,r}{\sum} \mathbb{E}[\Z_{1,k}^2\Z_{1,m}\Z_{2,m}\Z_{3,r} \Z_{4,r} ] \\
&+ N^2 \underset{m,r}{\sum} \mathbb{E}[\Z_{1,m}\Z_{2,m}\Z_{3,r} \Z_{4,r} ].
\label{eq:A1_A2_sum}
\end{align}
\\
 We recall that since $\Z_{i,k}$ and $\Z_{j,l}$ are independent for any $i$ and $j$ when $k \neq l$, we will always consider separately the cases where $k=l$ from the cases $k\neq l$. Let us first focus on the last term of \eqref{eq:A1_A2_sum}.

\begin{align*}
 \underset{m,r}{\sum} \mathbb{E}[\Z_{1,m}\Z_{2,m}\Z_{3,r} \Z_{4,r} ] & = \overset{N}{\underset{m=1}{\sum}} \mathbb{E}[\Z_{1,m}\Z_{2,m}\Z_{3,m}\Z_{4,m}] +  \underset{m \neq r}{\sum} \mathbb{E}[\Z_{1,m}\Z_{2,m}]\mathbb{E}[\Z_{3,r} \Z_{4,r} ] \\
 & = N \times o\left(\frac{1}{n}\right) + N(N-1) \times \frac{1}{(n-1)^2}
\end{align*}
\\
Then,

$\frac{1}{n^2} \frac{1}{N^4} \underset{i_1\neq i_2 \neq i_3 \neq i_4}{\sum} (N^2 \underset{m,r}{\sum} \mathbb{E}[\Z_{1,m}\Z_{2,m}\Z_{3,r} \Z_{4,r} ] ) = \frac{(N-1)^2(N-2)(N-3)}{n^2(n-1)^2} + o(1)$
\\
Now let us decompose the second term of \eqref{eq:A1_A2_sum} as:
\\
\begin{align*}
& \underset{k,m,r}{\sum} \mathbb{E}[\Z_{1,k}^2\Z_{1,m}\Z_{2,m}\Z_{3,r} \Z_{4,r} ] = \overset{N}{\underset{k=1}{\sum}} \mathbb{E}[\Z_{1,k}^3\Z_{2,k}\Z_{3,k} \Z_{4,k} ] + \underset{k \neq l}{\sum} \mathbb{E}[\Z_{1,k}^3\Z_{2,k}] \mathbb{E}[\Z_{3,l} \Z_{4,l} ]  \\
& +  \underset{k \neq l}{\sum} \mathbb{E}[\Z_{1,k}^2] \mathbb{E}[\Z_{1,l}\Z_{2,l}\Z_{3,l} \Z_{4,l} ] +  \underset{k \neq l}{\sum} \mathbb{E}[\Z_{1,k}^2\Z_{3,k}\Z_{3,k}] \mathbb{E}[\Z_{1,l} \Z_{2,l} ] + \underset{k \neq l \neq m}{\sum} \mathbb{E}[\Z_{1,k}^2] \mathbb{E}[\Z_{1,l}\Z_{2,l}] \mathbb{E}[\Z_{3,m} \Z_{4,m} ].
\end{align*}
\\
Using the results given by Proposition \ref{prop:prop_Z_EJS}, we obtain that 

$$\frac{1}{n^2}\frac{1}{N^4}\left(-2N \underset{k,m,r}{\sum} \mathbb{E}[\Z_{1,k}^2\Z_{1,m}\Z_{2,m}\Z_{3,r} \Z_{4,r} ] \right) = -\frac{2(N-1)^2(N-2)(N-3)}{n^2(n-1)^2} + o(1) .$$
\\
Similarly, we can prove that 

$$ \frac{1}{n^2}\frac{1}{N^4}\left( \underset{k,l,m,r}{\sum} \mathbb{E}[\Z_{1,k}^2 \Z_{2,l}^2\Z_{1,m}\Z_{2,m}\Z_{3,r} \Z_{4,r} ] \right) = \frac{(N-1)^2(N-2)(N-3)}{n^2(n-1)^2} + o(1), $$
\\
by using the properties of Proposition  \ref{prop:prop_Z_EJS} or similar relationships coming from other properties of $\Z$ that we have not detailed here.
\\
Hence we have shown \eqref{eq:A1_A2}.
The proofs of  \eqref{eq:A1_2}, \eqref{eq:A1_C12}, \eqref{eq:A1_C13}, \eqref{eq:C12_C23}, \eqref{eq:C13_C24}, \eqref{eq:C12_2} are very similar to this proof.
\\
It remains to prove Proposition \ref{prop3_card4}.
\\
According to the expression of $f_2(w,x,y,z)$ given in \eqref{eq:f2} and since

\begin{align*}
\vert \Sigma \vert ^{-\frac{1}{2}} = 1 -  \frac{\eta^\star}{2 \sqrt{N}}(A_1+A_2+A_3+A_4) & + \frac{\eta^{\star 2}}{4N} (A_1 A_2 + … + A_3 A_4) + \frac{3 \eta^{\star2}}{8N} (A_1^2+A_2^2+A_3^2+A_4^2) \\
& + \frac{\eta^{\star2}}{2N}(C_{1,2}^2+…+C_{3,4}^2)
+O_p\left(\frac{1}{N^\frac{3}{2}}\right),
\end{align*}
\\
the only term in $C_{1,2}C_{3,4}$ of $\mathbb{P}(\Y_{i_1}=\Y_{i_2}=\Y_{i_3}=\Y_{i_4}=1\vert \Z)$ is 
\\
$$\frac{1}{(2\pi)^2} \frac{\eta^{\star 2}}{N} \int_t^{+\infty}\int_t^{+\infty}\int_t^{+\infty}\int_t^{+\infty} wxyz C_{1,2}C_{3,4} dwdxdydz=\phi(t)^4  C_{1,2}C_{3,4} \frac{\eta^{\star 2}}{N}.$$
\\
The term in $C_{1,2}C_{3,4}$ of $\mathbb{P}(\textrm{``3 cases, 1 control"} \vert \Z)$ is
\\
$$-4\phi(t)^4 C_{1,2}C_{3,4}\frac{\eta^{\star 2}}{N}.$$
\\
The term in $C_{1,2}C_{3,4}$ of $\mathbb{P}(\textrm{``2 cases, 2 controls"} \vert \Z)$ is
\\
$$6 \phi(t)^4 C_{1,2}C_{3,4}\frac{\eta^{\star 2}}{N}.$$
\\
The term in $C_{1,2}C_{3,4}$ of $\mathbb{P}(\textrm{``1 case, 3 controls"} \vert \Z)$ is
\\
$$-4 \phi(t)^4 C_{1,2}C_{3,4}\frac{\eta^{\star 2}}{N}.$$
\\
The term in $C_{1,2}C_{3,4}$ of $\mathbb{P}(\Y_{i_1}=\Y_{i_2}=\Y_{i_3}=\Y_{i_4}=0\vert \Z)$ is
\\
$$ \phi(t)^4 C_{1,2}C_{3,4}\frac{\eta^{\star 2}}{N}.$$
\\
It remains to compute the approximation of the denominator of $\mathbb{E}[\W_{i_1} \W_{i_2} \W_{i_3} \W_{i_4}\vert \Z] $ of order 0, that is 

\begin{align*}
K^4 & +  4K^3(1-K) \left(\frac{K(1-P)}{P(1-K)}\right) + 6K^2 (1-K)^2  \left(\frac{K(1-P)}{P(1-K)}\right)^2  \\
& + 4K(1-K)^3 \left(\frac{K(1-P)}{P(1-K)}\right)^3 + (1-K)^4  \left(\frac{K(1-P)}{P(1-K)}\right)^4 \\
&= \frac{K^4}{P^4}\left[P^4 + 4 P^3(1-P) + 6P^2(1-P)^2+ 4P(1-P)^3+(1-P)^4\right] \\
&=\frac{K^4}{P^4}.
\end{align*}
\\
Finally, the term $C_{1,2}C_{3,4}$ in $\mathbb{E}[\W_{i_1} \W_{i_2} \W_{i_3} \W_{i_4}\vert \Z] $ is

\begin{align*}
\phi(t)^4 \frac{\eta^{\star2}}{N}C_{1,2}C_{3,4}& \left[ \frac{(1-P)^2}{P^2} + 2 \frac{1-P}{P} \left(\frac{K(1-P)}{P(1-K)}\right) + 6 \left(\frac{K(1-P)}{P(1-K)}\right)^2 \right.\\
& \left. + 2\frac{P}{1-P} \left(\frac{K(1-P)}{P(1-K)}\right)^3 + \left(\frac{K(1-P)}{P(1-K)}\right)^4 \right] \times \frac{P^4}{K^4} \\
&=\frac{P^2(1-P)^2}{K^4(1-K)^4} \phi(t)^4 \frac{\eta^{\star2}}{N}C_{1,2}C_{3,4},
\end{align*}
\\
which is exactly the term in $C_{1,2}C_{3,4}$ of  $\mathbb{E}[\W_{i_1} \W_{i_2}\vert \Z]\mathbb{E}[\W_{i_3} \W_{i_4}\vert \Z].$
\\
This proves Proposition \ref{prop3_card4}.


%
%
%
 
 \subsection{Second order approximation of $\mathbb{E}[\W_i \W_j \vert \Z, \epsilon_i=\epsilon_j=1]$}
 \label{pr:approx2}
 The density function $f$ can still be written as
 \begin{align*}
f(x,y) &= \frac{1}{2\pi \vert \Sigma^{(N)} \vert^{-\frac{1}{2}}} \exp \left\{ -\frac{1}{2 \vert \Sigma^{(N)} \vert}\left[x^2(1+\frac{\eta^\star}{\sqrt{N}}A_N(j)) +
y^2(1+\frac{\eta^\star}{\sqrt{N}}A_N(i)) -2xy \frac{B_N(i,j)}{\sqrt{N}} \right] \right\},
\end{align*}
  but with the explicit term of order $1/N$ in the expressions of $\vert \Sigma^{(N)} \vert^{-1}$ and $\vert \Sigma^{(N)} \vert^{-\frac{1}{2}}$:
  
   \begin{align*}
  \vert \Sigma^{(N)} \vert^{-1}& = 1 - \frac{\eta^\star}{\sqrt{N}}(A_N(i)+A_N(j)) + \frac{\eta^{\star 2}}{N}\left(-A_N(i) A_N(j) + B_N(i,j)^2 +(A_N(i)+A_N(j))^2\right)+O_p\left(\frac{1}{N^{\frac{3}{2}}}\right)\\
  &=1 - \frac{\eta^\star}{\sqrt{N}}(A_N(i)+A_N(j)) + \frac{\eta^{\star 2}}{N}\left(A_N(i) A_N(j) +A_N(i)^2+A_N(j)^2 + B_N(i,j)^2 \right)+O_p\left(\frac{1}{N^{\frac{3}{2}}}\right)
  \end{align*}
 and  
    \begin{align*}
\vert \Sigma^{(N)} \vert^{-\frac{1}{2}}& = 1 - \frac{\eta^\star}{2\sqrt{N}}(A_N(i)+A_N(j)) + \frac{\eta^{\star 2}}{2 N}\left(-A_N(i) A_N(j) + B_N(i,j)^2 +\frac{3}{8}(A_N(i)+A_N(j))^2\right) \\
& +O_p\left(\frac{1}{N^{\frac{3}{2}}}\right).
  \end{align*}
 \\ 
Thus, 

\begin{align*}
\exp & \left\{ -\frac{1}{2 \vert \Sigma^{(N)} \vert} \left[x^2(1+\frac{\eta^\star}{\sqrt{N}}A_N(j))  + y^2(1+\frac{\eta^\star}{\sqrt{N}}A_N(i)) -2 xy \frac{B_N(i,j)}{\sqrt{N}} \right] \right\}  \\
 & = \phi(x) \phi(y) \exp \left\{ -\frac{x^2}{2}(-A_N(i)\frac{\eta^\star}{\sqrt{N}} + \frac{\eta^{\star 2}}{N}(A_N(i)^2+B_N(i,j)^2)  \right. \\
 & \left. -\frac{y^2}{2}(-A_N(j)\frac{\eta^\star}{\sqrt{N}} + \frac{\eta^{\star 2}}{N}(A_N(j)^2+B_N(i,j)^2) ) +  xy (\frac{\eta^\star}{\sqrt{N}} B_N(i,j) - \frac{\eta^{\star 2}}{N}B_N(i,j)(A_N(i)+A_N(j))) \right\} \\
 & +O_p\left(\frac{1}{N^{\frac{3}{2}}}\right) \\
 & = \phi(x) \phi(y) \left[ 1 +  \frac{x^2}{2}(A_N(i)\frac{\eta^\star}{\sqrt{N}} - \frac{\eta^{\star 2}}{N}(A_N(i)^2+B_N(i,j)^2) +\frac{y^2}{2}(A_N(j)\frac{\eta^\star}{\sqrt{N}}\right. \\
 & - \frac{\eta^{\star 2}}{N}(A_N(j)^2+B_N(i,j)^2) )+ \frac{x^4}{8} \frac{\eta^{\star 2}}{N} A_N(i)^2 +\frac{y^4}{8} \frac{\eta^{\star 2}}{N} A_N(j)^2 \\
 & \left. + xy (\frac{\eta^\star}{\sqrt{N}} B_N(i,j) - \frac{\eta^{\star 2}}{N}B_N(i,j)(A_N(i)+A_N(j))) + \frac{x^2y^2}{2}\frac{\eta^{\star 2}}{N} B_N(i,j)^2+O_p\left(\frac{1}{N^{\frac{3}{2}}}\right) \right]
\end{align*}
\\
with the last term obtained by developing the exponential function.
\\
Since

\begin{align*}
 \int_{t}^{\infty} \int_{t}^{\infty} x^4 dx dy = t^3 \phi(t) + 3t\phi(t) + 3K
 \end{align*}
and

\begin{align*}
 \int_{t}^{\infty} \int_{t}^{\infty} x^2 y^2 dx dy =(t\phi(t) + K)^2,
 \end{align*}
 \\
 we have:
\begin{align*}
& \int_{t}^{\infty} \int_{t}^{\infty} \exp  \left\{ -\frac{1}{2 \vert \Sigma^{(N)} \vert} \left[x^2(1+\frac{\eta^\star}{\sqrt{N}}A_N(j))  + y^2(1+\frac{\eta^\star}{\sqrt{N}}A_N(i)) -2 xy \frac{B_N(i,j)}{\sqrt{N}} \right] \right\} dx dy \\
 &=  K^2 + \frac{K}{2}(t\phi(t)+K)\left[\frac{\eta^\star}{\sqrt{N}}(A_N(i)+A_N(j)) -\frac{\eta^{\star 2}}{N}(A_N(i)^2+A_N(j)^2+2 B_N(i,j)^2)\right] \\
 & + \frac{K}{8} \frac{\eta^{\star 2}}{N}(t^3\phi(t) +3t \phi(t)+3K) (A_N(i)^2+A_N(j)^2) + \phi(t)^2 \left[ \frac{\eta^{\star }}{\sqrt N} B_N(i,j)  - \frac{\eta^{\star 2}}{N} B_N(i,j)(A_N(i)+A_N(j))\right] \\
 &+ \frac{1}{2} (t\phi(t)+K)^2 \frac{\eta^{\star 2}}{N} (B_N(i,j)^2 +\frac{A_N(i) A_N(j)}{2}+  \frac{\phi(t)^2}{2} \frac{\eta^{\star 2}}{N} (t^2+2)B_N(i,j)(A_N(i)+A_N(j)) + O_p\left(\frac{1}{N^{\frac{3}{2}}}\right) \\
 \end{align*}
 \\
 Multiplying by

\begin{align*}
 \vert \Sigma^{(N)} \vert^{-\frac{1}{2}} & = 1 - \frac{\eta^\star}{2\sqrt{N}}(A_N(i)+A_N(j)) + \frac{\eta^{\star 2}}{2 N}\left(-A_N(i) A_N(j) + B_N(i,j)^2 \right. \\
 & \left. +\frac{3}{4}(A_N(i)+A_N(j))^2\right)+O_p\left(\frac{1}{N^{\frac{3}{2}}}\right), 
 \end{align*} 
\\ 
 we obtain 
 \\
 \begin{align*}
& \int_{t}^{\infty} \int_{t}^{\infty} f(x,y) dx dy 
 = K^2  + \frac{K}{2}t\phi(t)\frac{\eta^\star}{\sqrt{N}}(A_N(i)+A_N(j)) + \phi(t)^2  \frac{\eta^{\star }}{\sqrt N} B_N(i,j) + \frac{K}{8} \frac{\eta^{\star 2}}{N}(t^3\phi(t) -3 t \phi(t)) \\
 & + \frac{\eta^{\star 2}}{N} \frac{t^2\phi(t)^2}{4}A_N(i) A_N(j)+ \frac{\eta^{\star 2}}{N} B_N(i,j)^2 \frac{t^2}{2}\phi(t)^2 +\frac{\eta^{\star 2}}{N}\frac{\phi(t)^2}{2}(t^2-1) B_N(i,j) (A_N(i)+A_N(j)) \phi(t)^2 \\
 &+O_p\left(\frac{1}{N^{\frac{3}{2}}}\right).
  \end{align*}

 Similarly, 
 
\begin{align*}
 \int_{-\infty}^{t} \int_{-\infty}^{t} x^4 dx dy = -t^3 \phi(t) - 3t\phi(t) + 3(1-K)
 \end{align*}
and

\begin{align*}
 \int_{-\infty}^{t} \int_{-\infty}^{t} x^2 y^2 dx dy =(-t\phi(t) +1- K)^2.
 \end{align*}
 
  \begin{align*}
&  \int_{-\infty}^{t} \int_{-\infty}^{t} \exp  \left\{ -\frac{1}{2 \vert \Sigma^{(N)} \vert} \left[x^2(1+\frac{\eta^\star}{\sqrt{N}}A_N(j))  + y^2(1+\frac{\eta^\star}{\sqrt{N}}A_N(i)) -2 xy \frac{B_N(i,j)}{\sqrt{N}} \right] \right\} dx dy \\
&=(1-K)^2 +  \frac{1-K}{2}(-t\phi(t)+1-K)\left[\frac{\eta^\star}{\sqrt{N}}(A_N(i)+A_N(j)) -\frac{\eta^{\star 2}}{N}(A_N(i)^2+A_N(j)^2+2 B_N(i,j)^2)\right] \\
&+ \frac{1-K}{8} \frac{\eta^{\star 2}}{N}\left(-t^3\phi(t) -3t \phi(t)+3(1-K)\right) (A_N(i)^2+A_N(j)^2) + \frac{\phi(t)^2}{2} \frac{\eta^{\star 2}}{N} (t^2+2)B_N(i,j)(A_N(i)+A_N(j))
  \end{align*}
  
  $$ $$
  Multiplying by 
 
 $\vert \Sigma^{(N)} \vert^{-\frac{1}{2}} = 1 - \frac{\eta^\star}{\sqrt{N}}(A_N(i)+A_N(j)) + \frac{\eta^{\star 2}}{2 N}\left(-A_N(i) A_N(j) + B_N(i,j)^2 +\frac{3}{4}(A_N(i)+A_N(j))^2\right)+O_p\left(\frac{1}{N^{\frac{3}{2}}}\right)$,
 
  \begin{align*}
&  \int_{-\infty}^{t} \int_{-\infty}^{t} f(x,y) dx dy 
  =(1-K)^2 - \frac{1-K}{2}t\phi(t)\frac{\eta^\star}{\sqrt{N}}(A_N(i)+A_N(j)) + \phi(t)^2  \frac{\eta^{\star }}{\sqrt N} B_N(i,j) \\
  & + \frac{1-K}{8} \frac{\eta^{\star 2}}{N}(A_N(i)^2+A_N(j)^2)(-t^3\phi(t) +3 t \phi(t))  + \frac{\eta^{\star 2}}{N} \frac{t^2 \phi(t)^2}{4} A_N(i) A_N(j) \\
  & + \frac{\eta^{\star 2}}{N} B_N(i,j)^2 \frac{t^2}{2}\phi(t)^2 - \frac{\eta^{\star 2}}{N}B_N(i,j) (A_N(i)+A_N(j))\frac{ \phi(t)^2}{2} (t^2-1)+O_p\left(\frac{1}{N^{\frac{3}{2}}}\right).
  \end{align*}
  
  Finally, we compute similarly $\mathbb{P}(\Y_i \neq \Y_j \vert \Z)= \int_{-\infty}^{t} \int_{t}^{+\infty} f(x,y) dx dy +  \int_{t}^{+\infty} \int_{-\infty}^{t} f(x,y) dx dy$.
  
  We obtain
  
   \begin{align*}
\mathbb{P}(\Y_i \neq \Y_j \vert \Z) &=  2K(1-K) - \frac{1-2K}{2}t\phi(t)\frac{\eta^\star}{\sqrt{N}}(A_N(i)+A_N(j)) -2 \phi(t)^2  \frac{\eta^{\star }}{\sqrt N} B_N(i,j) \\
&+ \frac{1-2K}{8} \frac{\eta^{\star 2}}{N}(A_N(i)^2+A_N(j)^2)(t^3\phi(t) -3 t \phi(t))  - \frac{\eta^{\star 2}}{N}\frac{t^2\phi(t)^2}{2}A_N(i) A_N(j) \\
& - \frac{\eta^{\star 2}}{N} B_N(i,j)^2 t^2\phi(t)^2 + \frac{\eta^{\star 2}}{N}B_N(i,j) (A_N(i)+A_N(j)) \phi(t)^2 (-t^2+1)+O_p\left(\frac{1}{N^{\frac{3}{2}}}\right).
  \end{align*}

We replace the expressions of $\mathbb{P}(\Y_i =\Y_j=1 \vert \Z) $, $\mathbb{P}(\Y_i =\Y_j=0 \vert \Z) $ and $\mathbb{P}(\Y_i \neq \Y_j \vert \Z)$ in the expression of $\mathbb{E}[\W_i \W_j \vert \Z, \epsilon_i=\epsilon_j=1]$. Since we already computed the terms of order $\frac{1}{\sqrt N}$ for the numerator, it only remains the terms of order $\frac{1}{N}$.

Eventually, we find that the numerator can be writen as :

   \begin{align*}
\frac{\eta^\star}{\sqrt N} \frac{1-P}{P(1-K)^2}\phi(t)^2 B_N(i,j) & + \frac{\eta^{\star 2}}{N} \frac{t^2\phi(t)^2}{4} A_N(i) A_N(j)  \frac{1-P}{P(1-K)^2} + \frac{\eta^{\star 2}}{2N} B_N(i,j)^2 \frac{1-P}{P(1-K)^2}t^2\phi(t)^2  \\
&+ \frac{\eta^{\star 2}}{N} \frac{\phi(t)^2}{2}\frac{1-P}{P(1-K)^2} (t^2-1) B_N(i,j)(A_N(i)+A_N(j))+O_p\left(\frac{1}{N^{\frac{3}{2}}}\right).
  \end{align*}
  
 Similarly, we compute the expression of the denominator (at order $\frac{1}{\sqrt{N}}$ since the main term of the numerator is of order $\frac{1}{\sqrt{N}}$). We obtain the following expression:
 
   \begin{align*}
\frac{K^2}{P^2} +\frac{\eta^\star}{\sqrt N}  \frac{t}{2} \phi(t) (A_N(i)+A_N(j)) \frac{K(P-K)}{P^2(1-K)} + \frac{\eta^\star}{\sqrt N}  \phi(t)^2 B_N(i,j)\frac{(P-K)^2}{P^2(1-K)^2}+O_p\left(\frac{1}{N}\right).
  \end{align*}

     \begin{align*}
\mathbb{E}[\W_i \W_j \vert & \Z, \epsilon_i=\epsilon_j=1]= \frac{P^2}{K^2} \left[ 1-\frac{\eta^\star}{\sqrt N}  \frac{t}{2} \phi(t) (A_N(i)+A_N(j)) \frac{(P-K)}{K(1-K)}-\frac{\eta^\star}{\sqrt N}  \phi(t)^2 B_N(i,j)\frac{(P-K)^2}{K^2(1-K)^2}\right] \\
& \times \left[ \frac{\eta^\star}{\sqrt N} \frac{1-P}{P(1-K)^2}\phi(t)^2 B_N(i,j) + \frac{\eta^{\star 2}}{N} \frac{t^2\phi(t)^2}{4} A_N(i) A_N(j)  \frac{1-P}{P(1-K)^2}  \right. \\
&+ \left. \frac{\eta^{\star 2}}{2N} B_N(i,j)^2 \frac{1-P}{P(1-K)^2}t^2\phi(t)^2 - \frac{\eta^{\star 2}}{N} \frac{\phi(t)^2}{2}\frac{1-P}{P(1-K)^2} (t^2-1) B_N(i,j)(A_N(i)+A_N(j)) \right] \\
&+O_p\left(\frac{1}{N^{\frac{3}{2}}}\right) \\
&=  \frac{\eta^\star}{\sqrt N} \frac{P(1-P)}{K^2(1-K)^2}\phi(t)^2B_N(i,j) +\frac{t^2}{4} \frac{\eta^{\star 2}}{N}A_N(i) A_N(j) \frac{P(1-P)}{K^2(1-K)^2}  \\
&+  \frac{\eta^{\star 2}}{N} \frac{P(1-P)}{K^2(1-K)^2}\phi(t)^2B_N(i,j)^2 \left[\frac{t^2}{2}-\frac{(P-K)^2}{K^2(1-K)^2}\right] \\
& + \frac{\eta^{\star 2}}{2 N} \frac{P(1-P)}{K^2(1-K)^2}\phi(t)^2B_N(i,j)(A_N(i)+A_N(j)) \left[t^2-1 - \frac{P-K}{K(1-K)}t\phi(t) \right].
  \end{align*}
  
  \section*{Acknowledgements}
The author is very grateful to Elisabeth Gassiat for her insightful comments and valuable discussions.

\appendix  
\section{Appendix}  
   \subsection{Proof of Equation \eqref{eq:p_control}}
\label{app3}
By definition, the probabilities $p_{case}$ and $p_{control}$ are linked to the variables $\epsilon_i$ as follows:

$$p_{case}=\mathbb{P}(\epsilon_i=1 \vert \Z, \Y_i=1)$$

and 

$$p_{control}=\mathbb{P}(\epsilon_i=1 \vert \Z, \Y_i=0).$$

%
%

The ratio of the two following equations:
 $$P=\mathbb{P}(\Y_i=1 \vert \epsilon_i=1)= \frac{\mathbb{P}(\Y_i=1,\epsilon_i=1)}{\mathbb{P}(\epsilon_i=1)}=\frac{\mathbb{P}(\Y_i=1,V_i=1)}{\mathbb{P}(\epsilon_i=1)}=\frac{K p_{case}}{\mathbb{P}(\epsilon_i=1)}$$

and $$1-P=\mathbb{P}(\Y_i=0 \vert \epsilon_i=1)= \frac{\mathbb{P}(\Y_i=0,\epsilon_i=1)}{\mathbb{P}(\epsilon_i=1)}=\frac{\mathbb{P}(\Y_i=0,U_i=1)}{\mathbb{P}(\epsilon_i=1)}=\frac{(1-K) p_{control}}{\mathbb{P}(\epsilon_i=1)},$$

with the full ascertainment assumption given by \eqref{eq:full_asc} prove equation \eqref{eq:p_control}.

 \subsection{Proof of Equation \eqref{E_Wi_Wj}}
 \label{app1}
 This equation was proved in \cite{golan2014}, we recall the proof here for the sake of completeness.
 \\
 Conditionally to the event $\left\{\epsilon_i=\epsilon_j=1 \right\}$, the variable $\W_i\W_j$ can take the following values:
\\
$\bullet$ $\frac{1-p}{p}$ if $\Y_i=\Y_j=1$.
\\ 
$\bullet$  $\frac{p}{1-p}$ if $\Y_i=\Y_j=0$.
\\  
$\bullet$  $-1$ if $\Y_i \neq \Y_j$.

\vspace{5mm}

Let us write the expectaction of $\W_i \W_j$ conditionally to $\Z$ and conditionally to $\{\epsilon_i=\epsilon_j=1 \}$:

\begin{align}
\mathbb{E}[\W_i \W_j \vert \Z, \epsilon_i=\epsilon_j=1] & =  \frac{1-P}{P}  \mathbb{P}(\Y_i=\Y_j=1 \vert \Z, \epsilon_i=\epsilon_j=1)  - \mathbb{P}(\Y_i \neq \Y_j \vert \Z, \epsilon_i=\epsilon_j=1) \\ & + \frac{P}{1-P} \mathbb{P}(\Y_i=\Y_j=0 \vert \Z, \epsilon_i=\epsilon_j=1).
\label{eq:E1}
\end{align}

\begin{align*}
\mathbb{P}(\Y_i=\Y_j=1 \vert \Z, \epsilon_i=\epsilon_j=1) & = \frac{\mathbb{P}(  \epsilon_i=\epsilon_j=1\vert  \Y_i=\Y_j=1, \Z)\mathbb{P}(\Y_i=\Y_j=1 \vert \Z) }{\mathbb{P}(\epsilon_i=\epsilon_j=1 \vert \Z)} \\
& = \frac{\mathbb{P}(\Y_i=\Y_j=1 \vert \Z) }{\mathbb{P}(\epsilon_i=\epsilon_j=1 \vert \Z)}
\end{align*}
under the full ascertainment assumption given by Equation \eqref{eq:full_asc}.

Similarly, since we have seen in Equation \eqref{eq:p_control} that a control has a probability $\frac{K(1-P)}{P(1-K)}$ to be selected in the study and since $\epsilon_i$ and $\epsilon_j$ are assumed to be independent conditionally to $\Z$, $\Y_i$ and $\Y_j$: 
\begin{align*}
\mathbb{P}(\Y_i=\Y_j=0 \vert \Z, \epsilon_i=\epsilon_j=1) & = \frac{\mathbb{P}(  \epsilon_i=\epsilon_j=1\vert  \Y_i=\Y_j=0, \Z)\mathbb{P}(\Y_i=\Y_j=0 \vert \Z) }{\mathbb{P}(\epsilon_i=\epsilon_j=1 \vert \Z)} \\
& = \left(\frac{K(1-P)}{P(1-K)}\right)^2\frac{\mathbb{P}(\Y_i=\Y_j=1 \vert \Z) }{\mathbb{P}(\epsilon_i=\epsilon_j=1 \vert \Z)}
\end{align*}
and 
\begin{align*}
\mathbb{P}(\Y_i \neq \Y_j\vert \Z, \epsilon_i=\epsilon_j=1) & = \frac{\mathbb{P}(  \epsilon_i=\epsilon_j=1\vert  \Y_i \neq \Y_j, \Z)\mathbb{P}(\Y_i \neq \Y_j \vert \Z) }{\mathbb{P}(\epsilon_i=\epsilon_j=1 \vert \Z)} \\
& = \left(\frac{K(1-P)}{P(1-K)}\right)\frac{\mathbb{P}(\Y_i \neq \Y_j \vert \Z) }{\mathbb{P}(\epsilon_i=\epsilon_j=1 \vert \Z)}.
\end{align*}

The probability that both individuals $i$ and $j$ are included in the study is equal to
\begin{align*}
\mathbb{P}(\epsilon_i=\epsilon_j=1 \vert \Z) & =\mathbb{P}(\epsilon_i=\epsilon_j=1 \vert \Z, \Y_i=\Y_j=1) \mathbb{P}(\Y_i=\Y_j=1 \vert \Z) \\
& + \mathbb{P}(\epsilon_i=\epsilon_j=1 \vert \Z, \Y_i=\Y_j=0)\mathbb{P}(\Y_i=\Y_j=0 \vert \Z) + \mathbb{P}(\epsilon_i=\epsilon_j=1 \vert \Z, \Y_i\neq \Y_j)\mathbb{P}(\Y_i \neq \Y_j \vert \Z) \\
&=\mathbb{P}(\Y_i=\Y_j=1 \vert \Z) + \left(\frac{K(1-P)}{P(1-K)}\right)^2\mathbb{P}(\Y_i=\Y_j=0 \vert \Z) + \left(\frac{K(1-P)}{P(1-K)}\right)\mathbb{P}(\Y_i \neq \Y_j \vert \Z).
\end{align*}

If we combine all these computations and we plug them in the expression \eqref{eq:E1}, we obtain \eqref{E_Wi_Wj}.

\subsection{Proof of Equation \eqref{eq:Sigma_N}}
\label{app2}
Notice first that 
$$\G_N(i,i)-1= \frac{1}{\sqrt{N}}\left(\frac{1}{\sqrt{N}} \underset{k=1}{\overset{N}{\sum}} (\Z_{i,k}^2-1) \right)$$ with $$ \textrm{Var}\left(\frac{1}{\sqrt{N}} \underset{k=1}{\overset{N}{\sum}} (\Z_{i,k}^2-1) \right)= \frac{1}{N}\underset{k=1}{\overset{N}{\sum}} \mathbb{E}[\Z_{i,k}^4] - (\mathbb{E}[\Z_{i,k}^2])^2.$$

%
%

Moreover, since the variables $(\Z_{i,k})_{1\leq i \leq n}$ are normalized according to Equation \eqref{eq:normalization_1}, 

$$\overset{N}{\underset{i=1}{\sum}} \Z_{i,k}^2 =n.$$

By taking the expectation and since the variables $(\Z_{i,k})_{1\leq i \leq n}$ are exchangeable, we obtain that

\begin{equation}
\mathbb{E}[\Z_{i,k}^2]=1.
\label{eq:E_Zik}
\end{equation}

 Using \ref{prop1:2} of Proposition \ref{prop:prop_Z_EJS} and Equation \eqref{eq:E_Zik}, we obtain that $$ \textrm{Var}\left(\frac{1}{\sqrt{N}} \underset{k=1}{\overset{N}{\sum}} (\Z_{i,k}^2-1) \right)$$ is bounded and 
$$\G_N(i,i)-1=\frac{1}{N} \underset{k=1}{\overset{N}{\sum}} (\Z_{i,k}^2-1) = O_p\left(\frac{1}{\sqrt{N}}\right).$$ \\

Similarly,

\begin{align*}
 \textrm{Var}\left(\frac{1}{\sqrt{N}} \underset{k=1}{\overset{N}{\sum}} \Z_{i,k}\Z_{j,k} \right) & = \frac{1}{N}\underset{k=1}{\overset{N}{\sum}} \PE(\Z_{i,k}^2\Z_{j,k}^2)-\PE(\Z_{i,k}\Z_{j,k})^2 \\
 &= \frac{1}{N}\underset{k=1}{\overset{N}{\sum}} \left( 1 + o(1) - \frac{1}{(n-1)^2} \right) \textrm{ using \ref{prop1:3} and \ref{prop1:1} of Proposition \ref{prop:prop_Z_EJS}} \\
 & = 1+ o(1).
 \end{align*}
 
Then, $\frac{1}{N} \underset{k=1}{\overset{N}{\sum}} \Z_{i,k}\Z_{j,k} = O_p\left(\frac{1}{\sqrt{N}}\right)$.

Thus, we can write $$ \Sigma^{(N)}=\begin{pmatrix}
   1 + \frac{A_N(i)}{\sqrt{N}}\eta^\star &  \frac{B_N(i,j)}{\sqrt{N}}\eta^\star \\
    \frac{B_N(i,j)}{\sqrt{N}}\eta^\star & 1 + \frac{A_N(j)}{\sqrt{N}}\eta^\star
\end{pmatrix},
$$ 

 where $A_N(i)=\frac{1}{\sqrt{N}} \underset{k=1}{\overset{N}{\sum}} (\Z_{i,k}^2-1)=O_p(1)$ for all i,
 
  and  $B_N(i,j)=\frac{1}{\sqrt{N}} \underset{k=1}{\overset{N}{\sum}} \Z_{i,k}\Z_{j,k}=O_p(1)$ for all $i \neq j$.

%
%
%
%
%
 
 \subsection{Proof of Proposition \ref{prop:prop_Z_EJS}}
\label{app4}

Observe that for all $k=1,\dots,N$,
\begin{equation}\label{eq:sum_Z}
\sum_{i=1}^{n}\Z_{i,k}=0
\end{equation}
and
\begin{equation}\label{eq:sum_Z_car}
\sum_{i=1}^{n}\Z_{i,k}^{2}=n.
\end{equation}
Moreover, for each $k$, the random variables $(\Z_{i,k})_{1\leq i\leq n}$ are exchangeable. Thus, we deduce from
(\ref{eq:sum_Z_car}) that for all $i=1,\dots,n$ and $k=1,\dots,N$,
$
\PE(\Z_{i,k}^2)=1.
$
Hence, by (\ref{eq:sum_Z}), we get that
$$
0=\left(\sum_{i=1}^{n}\Z_{i,k}\right)^2=\sum_{i=1}^n \Z_{i,k}^2+\sum_{1\leq i\neq j\leq n} \Z_{i,k} \Z_{j,k}\;,
$$
which, by (\ref{eq:sum_Z_car}), implies that for all $k=1,\dots,N$ and $i\neq j=1,\dots,n$,
\begin{equation}\label{eq;prod_Z}
\PE(\Z_{i,k} \Z_{j,k})=-\frac{n}{n(n-1)}=-\frac{1}{n-1}\;,
\end{equation}
\\
that is \ref{prop1:1}.
\\
The proof of \ref{prop1:2} comes from the decomposition: 
\begin{align*}
\vert \Z_{1,k} \vert^p & =\vert \Z_{1,k} \vert^p \mathbbm{1}_{\{s_k^2 > \frac{\delta_{min}}{2} \}} + \vert \Z_{1,k} \vert^p \mathbbm{1}_{\{s_k^2 \leq \frac{\delta_{min}}{2} \}}\\
& \leq \frac{\vert A_{1,k} - \bar{A}_k \vert^p}{\left( \frac{\delta_{min}}{2}\right)^p} + n^p \mathbbm{1}_{\{s_k^2 \leq \frac{\delta_{min}}{2} \}}
\end{align*}
\\
Assumption \ref{hyp2} implies that $\underset{k}{\textrm{sup }} \mathbb{E}\left[ \vert A_{1,k} - \bar{A}_k \vert^p\right] < + \infty$ and 
the upper bound for $\mathbb{P}(s_k^2 \leq \delta)$ of Equation \eqref{eq:control_sk} prove \ref{prop1:2}.
\\
By (\ref{eq:sum_Z_car}), for all $k=1,\dots,N$, 
$$
n^2=\left(\sum_{i=1}^{n}\Z_{i,k}^{2}\right)^2=\sum_{i=1}^n \Z_{i,k}^{4}+\sum_{1\leq i\neq j\leq n} \Z_{i,k}^2 \Z_{j,k}^2\;
$$
Since the $(\Z_{i,k})_{1\leq i\leq n}$ are exchangeable for each $k=1,\dots,N$, we get that for all $k=1,\dots,N$,
$$
n=\PE[\Z_{1,k}^4]+(n-1)\PE[\Z_{1,k}^2 \Z_{2,k}^2]\;,
$$
which gives us \ref{prop1:3} by using \ref{prop1:2}.
\\
If we take the expectation of
$$\Z_{1,k}^3 \overset{n}{\underset{i=1}{\sum}} \Z_{i,k} =0,$$
we obtain 
$$\PE[\Z_{1,k}^4] +(n-1) \PE [\Z_{1,k}^3\Z_{2,k}] =0.$$

Then, \ref{prop1:2} implies \ref{prop1:4}.
\\
Similarly, since
$$\Z_{1,k} \Z_{2,k} \overset{n}{\underset{i=1}{\sum}} \Z_{i,k}^2 = n\Z_{1,k} \Z_{2,k}, $$
we obtain that 

$$2 \PE[\Z_{1,k}^3 \Z_{2,k}] + (n-2) \PE[\Z_{1,k}^2\Z_{2,k}\Z_{3,k}]= n \PE[\Z_{1,k} \Z_{2,k}].$$
\\
Then \ref{prop1:1} and  \ref{prop1:4} imply \ref{prop1:5}.
\\
Since 
$$ \Z_{1,k} \Z_{2,k} \Z_{3,k} \overset{n}{\underset{i=1}{\sum}} \Z_{i,k} =0,$$
we obtain that 

$$ 3 \PE [ \Z_{1,k}^2 \Z_{2,k} \Z_{3,k}] + (n-3) \PE[ \Z_{1,k} \Z_{2,k} \Z_{3,k} \Z_{4,k}]=0.$$
\\
Then, \ref{prop1:5} implies \ref{prop1:6}. 
\\
Since 
$$ \Z_{1,k}^5 \overset{n}{\underset{i=1}{\sum}} \Z_{i,k} =0,$$
we obtain that 

$$  \PE [ \Z_{1,k}^6] + (n-1) \PE[ \Z_{1,k}^5 \Z_{2,k} ]=0.$$
Then, \ref{prop1:2} implies \ref{prop1:7}. 
\\
The proof of \ref{prop1:8} is very similar to the proof of \ref{prop1:2} but we use Assumption \ref{hyp4} which gives us that 
$\underset{k}{\textrm{sup }} \mathbb{E}\left[ \vert (A_{1,k} - \bar{A}_k)(A_{2,k} - \bar{A}_k) \vert^p\right] < + \infty$.
\\
 Since $$ \Z_{1,k}^4 \overset{n}{\underset{i=1}{\sum}} \Z_{i,k}^2= n \Z_{1,k}^4,$$
 
 $$ \mathbb{E}[\Z_{1,k}^6 ]+ (n-1) \mathbb{E}[\Z_{1,k}^4\Z_{2,k}^2]=n \mathbb{E}[\Z_{1,k}^4].$$
 \\
 Then \ref{prop1:2} implies \ref{prop1:9}.
 \\
 Similarly, since 
 $$\Z_{1,k}^4 \Z_{2,k} \overset{n}{\underset{i=1}{\sum}} \Z_{i,k} =0,$$ 
 we obtain that

$$ \mathbb{E}[\Z_{1,k}^5\Z_{2,k}] + \mathbb{E}[\Z_{1,k}^4\Z_{2,k}^2] + (n-2) \mathbb{E}[\Z_{1,k}^4\Z_{2,k}\Z_{3,k}] =0.$$
Then, \ref{prop1:7} and \ref{prop1:9} imply \ref{prop1:10}.
\\
Since
 $$\Z_{1,k}^3 \Z_{2,k} \overset{n}{\underset{i=1}{\sum}} \Z_{i,k}^2 = n\Z_{1,k}^3 \Z_{2,k},$$
 we obtain that 
 
 $$\mathbb{E}[\Z_{1,k}^5\Z_{2,k}] +  \mathbb{E}[\Z_{1,k}^3\Z_{2,k}^3] + (n-2) \mathbb{E}[\Z_{1,k}^3\Z_{2,k}^2\Z_{3,k}] = n \mathbb{E}[\Z_{1,k}^3\Z_{2,k}].$$
 Then, \ref{prop1:7}, \ref{prop1:8} and \ref{prop1:4} imply \ref{prop1:11}.
\\
 Finally, since $\Z_{1,k}^3 \Z_{2,k} (\overset{n}{\underset{i=1}{\sum}} \Z_{i,k})^2=0$,

\begin{align*}
& \mathbb{E}[\Z_{1,k}^5\Z_{2,k}]  +  \mathbb{E}[\Z_{1,k}^3\Z_{2,k}^3]  + 2 \mathbb{E}[\Z_{1,k}^4\Z_{2,k}^2] + 2(n-2)  \mathbb{E}[\Z_{1,k}^4\Z_{2,k}\Z_{3,k}]  \\
& + 2(n-2) \mathbb{E}[\Z_{1,k}^3\Z_{2,k}^2\Z_{3,k}]  +(n-2)^2 \mathbb{E}[\Z_{1,k}^3\Z_{2,k}\Z_{3,k}\Z_{4,k}] =0
\end{align*}
Then, \ref{prop1:7},  \ref{prop1:8}, \ref{prop1:9}, \ref{prop1:10} and \ref{prop1:11} imply \ref{prop1:12}.

  \bibliographystyle{abbrvnat-mod}
  \bibliography{bibli_intro.bib}
  
\end{document}